\documentclass[11pt]{amsart}
\numberwithin{equation}{section}
\newtheorem{theorem}{Theorem}

\newtheorem{lemma}{Lemma}
\newtheorem{corollary}{Corollary}
\newtheorem{remark}{Remark}[section]

\numberwithin{theorem}{section} \numberwithin{lemma}{section}
\numberwithin{proposition}{section} \numberwithin{equation}{section}

\begin{document}

\tracingpages 1
\title[gradient estimate]{\bf A gradient estimate for all positive solutions
of the conjugate heat equation under Ricci flow}
\author{ Shilong Kuang and Qi S. Zhang}
\address{Department of Mathematics,  University of California,
Riverside, CA 92521, USA }
\date{October 2006}
\begin{abstract}
We establish a point-wise gradient estimate for $all$ positive
solutions of the conjugate heat equation. This contrasts to
Perelman's point-wise gradient estimate which works mainly for the
fundamental solution rather than all solutions. Like Perelman's
estimate, the most general form of our gradient estimate does not
require any curvature assumption. Moreover, assuming only lower
bound on the Ricci curvature, we also prove a localized gradient
estimate similar to the Li-Yau estimate for the linear
Schr\"odinger heat equation. The main difference with the linear
case is that no assumptions on the derivatives of the potential
(scalar curvature) are needed.

A generalization of Perelman's W-entropy is defined in both the
Ricci flow and fixed metric case. We also find a new family of heat
kernel estimates.
\end{abstract}

\maketitle
\tableofcontents
\section{Introduction}


In the paper [P], Perelman discovers a monotonicity formula for
the $W$ entropy of positive solutions of the conjugate heat
equation.
\begin{equation}
\label{coneq}
\begin{cases}
\Delta u - R u + \partial_t u=0\\
\partial_t g = - 2 Ric.
\end{cases}
\end{equation}
Here $u=u(x, t)$ with $x \in {\bf M}$, a compact manifold and $t \in
(0, T)$, $T>0$;  $\Delta$ is the Laplace-Beltrami operator under the
metric $g$ and $Ric$ is the Ricci curvature tensor. Here and though
out it is assumed the metric $g$ is smooth in the region ${\bf M}
\times (0, T)$ with $T>0$, unless stated otherwise. Moreover, he
shows that this formula implies a point-wise gradient estimate for
the fundamental solution of the conjugate heat equation (Corollary
9.3 [P]). Namely, let $u$ be the fundamental solution of
(\ref{coneq}) in ${\bf M} \times (0, T)$ and $f$ be the function
such that $u=(4 \pi \tau)^{-n/2} e^{-f}$ with $\tau = T-t$. Then
\[
[\tau(2\triangle f-|\nabla f|^2+R)+f-n]u \le 0
\] in ${\bf M} \times (0, T)$.
This formula can be regarded as a generalization of the
Li-Yau-Hamilton gradient estimate for the heat equation. By now it
is clear that the importance of Perelman's monotonicity formula
and gradient estimate can hardly be overstated. See for example
\cite{CCGGIIKLLN:1}, \cite{CZ:1}, \cite{KL:1} and \cite{MT:1}.
However, there is one place where some improvement is still
desirable, namely the gradient estimate does not apply to all
positive solutions to the conjugate heat equation. For instance,
for the Ricci flat manifold $S^1 \times S^1$. The constant $1$ is
a solution to the conjugate heat equation. Clearly it does not
satisfy Perelman's gradient estimate stated above.  Whether a
Perelman type gradient estimate exists for all positive solutions
of the conjugate heat equation is a question circulating for a few
years.

The main goal of this paper is to establish a gradient estimate that
works for $all$ positive solutions of the conjugate heat equation.
Like Perelman's estimate for the fundamental solution, the most
general form of the new gradient estimate does not require any
curvature assumption. Moreover, assuming only lower bound on the
Ricci curvature, it also has a local version which appears similar
to the Li-Yau estimate for the linear heat equation. An immediate
consequence of the gradient estimate is a classical Harnack
inequality for positive solutions of the conjugate heat equation.

We also introduce a generalization of Perelman's $W$ entropy and
prove its monotonicity. Specializing to the fixed metric case, we
prove a family of gradient estimates for the fundamental solution
of the heat equation. This family includes the Li-Yau estimate and
Perelman's estimate (specialized to the heat equation cf
\cite{N1}) as special cases.

The rest of the paper is organized as follows.  The results
concerning the conjugate heat equation under Ricci flow is given in
sections 2-3. In section 4 we set up a generalization of Perelman's
W-entropy. In Section 5, we present the results for the linear heat
equation in the fixed metric case. Some useful calculations which
are around in various papers and preprints are collected in the
appendix.

\section{New Gradient Estimate and Harnack Inequality
 for positive solutions to the conjugate heat equation}


The main result of this section is

\begin{theorem}  \label{Li-Yau-Perelman}
Suppose g(t) evolve by the Ricci flow, that is, $\frac{\partial
g}{\partial t}=-2 Ric $ on a closed manifold $M$ for $t\in [0,T)$,
and $u: M\times[0,T) \mapsto(0,\infty)$ be a positive $C^{2,1}$
solution to the conjugate heat equation $\Box^{*}u=-\triangle
u-u_{t}+Ru=0$. Let $u=\frac{e^{-f}}{(4\pi\tau)^\frac{n}{2}}$
 and $\tau=T-t$. Then: \\
(i) if the scaler curvature $R\geq 0$, then for all $t \in (0, T)$
and all points,
\begin{equation}  \label{Li-Yau-Perelman-inequality}
 2\triangle f -|\nabla f|^2+R\leq \frac{2n}{\tau};
\end{equation}
\\
(ii) without assuming the non-negativity of $R$, then for
$t\in[\frac{T}{2},T)$ and all points,
\begin{equation}  \label{Li-Yau-Perelman-inequality-1}
 2\triangle f -|\nabla f|^2+R\leq \frac{3n}{\tau}.
\end{equation}

\end{theorem}
%
%
\begin{remark} Since $f=-\ln u-\frac{n}{2}\ln(4\pi\tau)$, if we
replace $f$ by $u$ accordingly, then  we get
\begin{equation}
\begin{aligned}
 &\frac{|\nabla u|^2}{u^2}
 -2\frac{u_{\tau}}{u}-R\leq \frac{2n}{\tau}, \quad \text{if} \quad R \ge
 0;
\\
 &\frac{|\nabla u|^2}{u^2}
 -2\frac{u_{\tau}}{u}-R\leq \frac{3n}{\tau}, \quad \text{if $R$
 changes sign and} \quad t \ge T/2.
\end{aligned}
\end{equation}
It is similar to the Li-Yau gradient estimate  for the heat
equation on manifolds with nonnegative Ricci curvature, i.e.
\[
\frac {|\nabla u|^2}{u^2}-\frac{u_t}{u}\leq\frac{n}{2t}
\]for positive solutions  of $\Delta u - \partial_t u =0.$
\end{remark}

\begin{remark} Some related gradient estimates with various
dependence on the Ricci and other curvatures can be found in
\cite{G} and \cite{N3}.
\end{remark}

\begin{proof} of Theorem 2.1.

By a standard approximation argument, we can assume without loss
of generality that $g=g(t)$ is smooth in the closed time interval
$[0, T]$ and that $u$ is strictly positive everywhere.

(i) By  standard computation (one can consult various sources for
more details(\cite{CK} and \cite{T} e.g.)),
\begin{align}
    \left(
    \frac{\partial}{\partial t}+\triangle
    \right)
    (\triangle f)
 &=\triangle \frac{\partial f}{\partial t}
   +2\langle Ric,Hess(f)\rangle
   +\triangle (\triangle f)
\notag\\
 &=\triangle
   \left(
   -\triangle f +|\nabla f|^2-R+\frac{n}{2\tau}
   \right)
   +2\langle Ric,Hess(f)\rangle
\notag\\
  &\quad+\triangle (\triangle f)
\notag\\
  &=2\langle Ric,Hess(f)\rangle +\triangle
   \left(
   |\nabla f|^2-R
   \right).
 \end{align}
Also using the evolution equation of $g$,
\begin{align}
    \left(
    \frac{\partial}{\partial t}+\triangle
    \right)
    |\nabla f|^2
 &=2Ric(\nabla f,\nabla f)
  +2\langle\,\nabla f,\nabla \frac{\partial f}{\partial t}\rangle
  +\triangle \, |\nabla f|^2
 \notag\\
 &=2Ric(\nabla f,\nabla f)
  +2\langle\,\nabla f,\nabla (-\triangle f +|\nabla f|^2-R)\rangle
\notag\\
  &\quad
   +\triangle \, |\nabla f|^2 \,.
\end{align}
Notice also
\begin{equation}
    \left(
    \frac{\partial}{\partial t}+\triangle
    \right)R
=2\triangle R +2|Ric|^2 \,.
\end{equation}
Combining these three expressions, we deduce
\begin{align}    \label{q}
    \left(
    \frac{\partial}{\partial t}+\triangle
    \right)
   & (2\triangle f-|\nabla f|^2+R)
\notag\\
&=4\langle\,Ric,Hess(f)\rangle +\triangle |\nabla f|^2-2Ric(\nabla
f,\nabla f)
\notag\\
 &\quad
  -2\langle\,\nabla f,\nabla(-\triangle f +|\nabla f|^2-R)\rangle
  +2|Ric|^2 \,.
\end{align}
Denote
 \[
 q(x,t)=2\triangle f-|\nabla f|^2+R.
 \]By Bochner's identity,
 \[
 \triangle
|\nabla f|^2=2|f_{ij}|^2+2\nabla f \nabla (\triangle
f)+2R_{ij}f_{i}f_{j},
\]
the above equation becomes
\begin{align}
    \left(
    \frac{\partial}{\partial t}+\triangle
    \right)
    q &=4 R_{ij}f_{ij}+\big(2|f_{ij}|^2+2\nabla f\, \nabla (\triangle
f)+2R_{ij}f_{i}f_{j}\big)-2R_{ij}f_{i}f_{j}
\notag\\
 &\quad
  -2\nabla f\,\nabla(-\triangle f +|\nabla f|^2-R)+2R_{ij}^2
\notag\\
 &=4 R_{ij}f_{ij}+2|f_{ij}|^2+2R_{ij}^2+2\nabla f \,\nabla \big(2\triangle
f-|\nabla f|^2+R\big) \notag\\
 &=2|R_{ij}+f_{ij}|^2+2\nabla f \,\nabla q
\notag
\end{align}
that is,
\begin{align}
    \left(
    \frac{\partial}{\partial t}+\triangle
    \right)
    q-2\nabla f \,\nabla q
    &=2|R_{ij}+f_{ij}|^2
    \geq \frac{2}{n}(R+\triangle f)^2  \,.
\end{align}
Since
\[
q=2\triangle f-|\nabla f|^2+R=2(\triangle f+R)-|\nabla
f|^2-R,
\]
and hence
\[
R+\triangle f=\frac{1}{2}(q+|\nabla f|^2+R),
\]
we have
\begin{align}  \label{heat-operator-q}
    \left(
    \frac{\partial}{\partial t}+\triangle
    \right)
    q-2\nabla f \,\nabla q
    \geq \frac{1}{2n}\big(q+|\nabla f|^2+R\big)^2 \,.
\end{align}
By direct computation, we also have, for any $\epsilon >0$
\begin{align}
    \left(
    \frac{\partial}{\partial t}+\triangle
    \right)
     \frac{2n}{T-t+\epsilon} -2\nabla f \,\nabla \big( \frac{2n}{T-t+\epsilon} \big)
    = \frac{1}{2n}\big( \frac{2n}{T-t+\epsilon} \big)^2 \,.
\end{align}
Combine the above two expressions, we get
\begin{align}    \label{maximum-principle-inequality}
 &\quad
    \left(
          \frac{\partial}{\partial t}+\triangle
    \right)
     \big(   q-\frac{2n}{T-t+\epsilon}\big) -2\nabla f \,\nabla \big(q- \frac{2n}{T-t+\epsilon}    \big)
\notag\\
 & \geq \frac{1}{2n}
     \left(
           q+\frac{2n}{T-t+\epsilon}+|\nabla f|^2+R
     \right)
     \left(
           q-\frac{2n}{T-t+\epsilon}+|\nabla f|^2+R
     \right) \,.
\end{align}
We deal with the above inequality in two cases:\\
\\
\textbf{Case 1.} at a point $(x, t)$,
$q+\frac{2n}{T-t+\epsilon}+|\nabla f|^2+R\leq0$, then also
\[
 q-\frac{2n}{T-t+\epsilon}+|\nabla f|^2+R \leq0
\]
thus,
\begin{align}
    &\quad\left(
    \frac{\partial}{\partial t}+\triangle
    \right)
     \big(q-\frac{2n}{T-t+\epsilon}\big) -2\nabla f \,\nabla \big(q- \frac{2n}{T-t+\epsilon} \big)
    \geq 0.
\end{align}
\\
\textbf{Case 2.} at a point $(x, t)$,
$q+\frac{2n}{T-t+\epsilon}+|\nabla f|^2+R>0$, then the inequality
(\ref{maximum-principle-inequality}) can be transformed to
\begin{align}
    &\quad
    \left(
        \frac{\partial}{\partial t}+\triangle
    \right)
     \big(  q-\frac{2n}{T-t+\epsilon}\big) -2\nabla f \,\nabla \big(q- \frac{2n}{T-t+\epsilon}   \big)
\notag\\
    & \quad- \frac{1}{2n}  \left(    q+\frac{2n}{T-t+\epsilon}+|\nabla f|^2+R    \right)
    \left(     q-\frac{2n}{T-t+\epsilon}    \right)
\notag\\
    & \quad
     \geq \frac{1}{2n}
     (|\nabla f|^2+R)
     \left(    q+\frac{2n}{T-t+\epsilon}+|\nabla f|^2+R    \right)
    \geq0 \,.
\end{align}
Defining a potential term by
\begin{equation}
 V=V(x, t)=\begin{cases}
    0;   &\text{if $q+\frac{2n}{T-t+\epsilon}+|\nabla f|^2+R \leq 0$
    at $(x, t)$}\\
    \frac{1}{2n}
    \big(
         q+\frac{2n}{T-t+\epsilon}+|\nabla f|^2+R
    \big)
    ; &\text{if $q+\frac{2n}{T-t+\epsilon}+|\nabla f|^2+R \geq 0$
    at $(x, t)$} \,.
  \end{cases}
 \end{equation}
We know $V$ is continuous; further, by the above two cases, we
conclude
\begin{align}
    \big(
        \frac{\partial}{\partial t}+\triangle
    \big)
     \big(  q-\frac{2n}{T-t+\epsilon}\big) -2\nabla f \nabla \big(q- \frac{2n}{T-t+\epsilon}   \big)
    -V
    \big(    q-\frac{2n}{T-t+\epsilon}   \big)
    \geq0.
\end{align}
Since we assumed that the Ricci flow is smooth in $[0,T]$ and that
$u(x,t)$ is a positive $C^{2,1}$ solution to the conjugate heat
equation, thus
\[
q=2\triangle f-|\nabla f|^2+R =\frac{|\nabla
u|^2}{u^2}-\frac{2\triangle u}{u}+R
\] is bounded for $t\in[0,T]$.
If we choose $\epsilon$ sufficiently small, then
$q(x,T)\leq\frac{2n}{\epsilon}$, thus by the maximum
principle(\cite{CK}, e.g.), for all $t\in[0,T]$,
$q(x,t)\leq\frac{2n}{T-t+\epsilon}$. Let $\epsilon \rightarrow 0$,
we have for all $t\in[0,T]$,
\begin{align}
q(x,t)\leq\frac{2n}{T-t} \,.
\end{align}
Recall $q=2\triangle f-|\nabla f|^2+R$, $\tau=T-t$, then we have
\begin{equation}
 2\triangle f -|\nabla f|^2+R\leq \frac{2n}{\tau} \,.
\end{equation}
Further, $f=-\ln u-\frac{n}{2}(4\pi\tau)$, then the above yields
\begin{equation}
 \frac{|\nabla u|^2}{u^2}-\frac{2u_{\tau}}{u}-R
 \leq
 \frac{2n}{\tau}.
\end{equation}
\\
\medskip

{\bf Proof of (ii).}  Next we prove the gradient estimate without
the non-negativity assumption for the scalar curvature $R$. Let
$c\geq 2n$ be a constant to be determined later; denote
\[
B=|\nabla f|^2+R  \,.
\]
Similar to the inequality (\ref{maximum-principle-inequality}), we
also have,
\begin{align} \label{maximum-principle-inequality-1}
 &\quad
    \left(
          \frac{\partial}{\partial t}+\triangle
    \right)
     \big(   q-\frac{c}{T-t+\epsilon}\big) -2\nabla f \,\nabla \big(q- \frac{c}{T-t+\epsilon}    \big)
\\
 &\quad \geq \frac{1}{2n}
     \left(
           q+B
     \right)^2
     -\frac{c}{(T-t+\epsilon)^2}
\notag\\
&\quad = \frac{1}{2n}
   [
     (q+B)^2
     -\frac{c^2}{(T-t+\epsilon)^2}
     +\frac{c^2}{(T-t+\epsilon)^2}
     -\frac{2cn}{(T-t+\epsilon)^2}
   ]
\notag\\
&\quad = \frac{1}{2n}
   [
     (q-\frac{c}{T-t+\epsilon}+B)
     (q+\frac{c}{T-t+\epsilon}+B)
     +\frac{c(c-2n)}{(T-t+\epsilon)^2}
   ].
\notag
\end{align}
We deal with the previous inequality at a given point $(x,t)$ in three cases :\\
\\
\textbf{Case 1.} $B\geq0$, and $q+\frac{c}{T-t+\epsilon}+B\leq0$,
then also
\[
 q-\frac{c}{T-t+\epsilon}+B \leq0
\]
thus,
\begin{align}
    &\quad\left(
    \frac{\partial}{\partial t}+\triangle
    \right)
     \big(q-\frac{c}{T-t+\epsilon}\big) -2\nabla f \,\nabla \big(q- \frac{c}{T-t+\epsilon} \big)
    \geq 0.
\end{align}
\\
\textbf{Case 2.}
$B\geq0$, and $q+\frac{c}{T-t+\epsilon}+B>0$, then
the inequality (\ref{maximum-principle-inequality-1}) can be changed
to
\begin{align}
    &\quad
    \left(
        \frac{\partial}{\partial t}+\triangle
    \right)
     \big(  q-\frac{c}{T-t+\epsilon}\big) -2\nabla f \,\nabla \big(q- \frac{c}{T-t+\epsilon}   \big)
\notag\\
    & \quad-\frac{1}{2n}  \left(    q+\frac{c}{T-t+\epsilon}+B   \right)
    \left(     q-\frac{c}{T-t+\epsilon}    \right)
\notag\\
    & \quad
     \geq \frac{1}{2n}
     B
     \left(    q+\frac{c}{T-t+\epsilon}+B    \right)
    \geq0.
\end{align}
%
\textbf{Case 3.}
$B\leq0$, then the inequality
(\ref{maximum-principle-inequality-1}) can be changed to
\begin{align}
    &\quad
    \left(
        \frac{\partial}{\partial t}+\triangle
    \right)
     \big(  q-\frac{c}{T-t+\epsilon}\big) -2\nabla f \,\nabla \big(q- \frac{c}{T-t+\epsilon}   \big)
\notag\\
    & \quad \geq
    \frac{1}{2n}  \left(    q+\frac{c}{T-t+\epsilon}+B   \right)
    \left(     q-\frac{c}{T-t+\epsilon}    \right)
\notag\\
    & \quad\quad
     + \frac{1}{2n}
     B
     \left(    q-\frac{c}{T-t+\epsilon}    \right)
     +\frac{1}{2n}
     \left(
           \frac{2Bc}{T-t+\epsilon}+\frac{c(c-2n)}{(T-t+\epsilon)^2}
     \right) \,.
\end{align}
To continue, we need the following estimate of the scalar
curvature $R$ under the Ricci flow i.e.
\begin{align}  \label{estimate of R}
R\geq -\frac{n}{2(t+\epsilon)}
\end{align}
for some $\epsilon>0$ depending on the initial value of $R$(\,it
comes from the weak minimum principle for a differential inequality
$\frac{\partial R}{\partial t}\geq \triangle R+ \frac{2}{n}R^2$
(\cite{CK}, e.g.); thus
\begin{align}
B=|\nabla f|^2+R\geq R\geq
-\frac{n}{2(t+\epsilon)}\geq-\frac{n}{2(T-t+\epsilon)}
\end{align}
for $t\geq\frac{T}{2}$ because
$t\geq\frac{T}{2}
 \,\,\Rightarrow\,\,
 t\geq T-t
 \,\,\Rightarrow\,\,
 t+\epsilon
 \geq
 T-t+\epsilon
 \,\,\Rightarrow
 \,\,\frac{1}{t+\epsilon}\leq\frac{1}{T-t+\epsilon}$
;\\
thus
\begin{align}
 &\frac{1}{2n}
     \left( \frac{2Bc}{T-t+\epsilon}+\frac{c(c-2n)}{(T-t+\epsilon)^2}     \right)
\notag\\
 &\geq
 \frac{1}{2n}
     \left(-\frac{n}{2(T-t+\epsilon)} \, \frac{2c}{T-t+\epsilon}+\frac{c(c-2n)}{(T-t+\epsilon)^2}     \right)
\notag\\
 &=
 \frac{1}{2n}
     \left(\frac{c(c-3n)}{(T-t+\epsilon)^2}     \right) \,.
\end{align}
Therefore
\begin{align}
    &\quad
    \left(
        \frac{\partial}{\partial t}+\triangle
    \right)
     \big(  q-\frac{c}{T-t+\epsilon}\big) -2\nabla f \,\nabla \big(q- \frac{c}{T-t+\epsilon}   \big)
\notag\\
    & \quad -
    \frac{1}{2n}  \left(    q+\frac{c}{T-t+\epsilon}+2B   \right)
    \left(     q-\frac{c}{T-t+\epsilon}    \right)
\notag\\
    &\quad\geq \frac{c\,(c-3n)}{2n(T-t+\epsilon)^2}
\end{align}
take $c=3n$, we have,
\begin{align}
    \big(
        \frac{\partial}{\partial t}+\triangle
    \big)
     \big(  q-\frac{2n}{T-t+\epsilon}\big) -2\nabla f \nabla \big(q- \frac{2n}{T-t+\epsilon}   \big)
    -V
    \big(    q-\frac{2n}{T-t+\epsilon}   \big)
    \geq0
\end{align}
where $V=V(x, t)$ is a continuous function defined by
\begin{equation}
 V=\begin{cases}
    0;   &\text{if $B\geq 0$, $q+\frac{2n}{T-t+\epsilon}+B \leq 0$ at
    $(x, t)$}\\
    \frac{1}{2n}
    \big(
         q+\frac{2n}{T-t+\epsilon}+B
    \big);
    &\text{if $B\geq 0$, $q+\frac{2n}{T-t+\epsilon}+B > 0$
    at
    $(x, t)$}\\
    \frac{1}{2n}
    \big(
         q+\frac{2n}{T-t+\epsilon}+2B
    \big);
    &\text{if $B< 0$ at
    $(x, t)$} \,.
    \end{cases}
 \end{equation}
Follow the similar argument for the inequality
(\ref{Li-Yau-Perelman-inequality}), by the maximum principle
again, we have
\begin{equation}
 2\triangle f -|\nabla f|^2+R\leq \frac{3n}{\tau}\quad and\quad \frac{|\nabla u|^2}{u^2}-\frac{2u_{\tau}}{u}-R
 \leq
 \frac{3n}{\tau}, \quad t \ge T/2.
\end{equation}
\end{proof}

\medskip
An immediate consequence of the above theorem is:

\begin{corollary}[Harnack Inequality]
Given a smooth Ricci flow on a closed manifold $M$, let
$u:M\times[0,T) \mapsto(0,\infty)$ be a positive $C^{2,1}$ solution
to the conjugate heat equation.

(a). Suppose the scalar curvature $R \ge 0$ for $t\in [0,T)$. Then
for any two points $(x, t_1)$, $(y, t_2)$ in ${\bf M} \times (0, T)$
such that $t_1 < t_2$, it holds
\begin{align}
u(y, t_2)
 \leq
 u(x, t_1)
 \left(
            \frac{\tau_{1}}{\tau_{2}}
       \right)^n
 \exp
  {
   \frac
        {\int_{0}^{1}      [\,4|\gamma^{\prime}(s)|^2+(\tau_1-\tau_2)^2\,R\,]\,    ds}
        {2(\tau_1-\tau_2)}
  }.
\end{align}
Here $\tau_i = T-t_i$, $i =1, 2$, and $\gamma(s):[0,1]\rightarrow M$
is a smooth curve from $x$ to $y$.

(b). Without assuming the nonnegativity of the scalar curvature $R$,
then for $t_2>t_1 \ge T/2$, it holds
\begin{align}
u(y, t_2)
 \leq
 u(x, t_1)
 \left(
            \frac{\tau_{1}}{\tau_{2}}
       \right)^{3n/2}
 \exp
  {
   \frac
        {\int_{0}^{1}      [\,4|\gamma^{\prime}(s)|^2+(\tau_1-\tau_2)^2\,R\,]\,    ds}
        {2(\tau_1-\tau_2)}
  }.
\end{align}

\end{corollary}

\begin{proof} We will only prove (a) since the proof of (b) is
similar.

Denote $\tau(s):=\tau_{2}+(1-s)(\tau_{1}-\tau_{2})$, $0\leq
\tau_{2}<\tau_{1}\leq T$, define
\[
\ell(s):=\ln \,u
                (\gamma(s),T-\tau(s))
\]
where $\ell(0)=\ln\,u(x, t_1)$, $\ell(1)=\ln\,u(y, t_2)$.
\\
By direct computation,
\begin{align}
\frac{\partial \ell(s)}{\partial s}
 &= \frac{u_{s}}{u}
  = \frac{\nabla u}{u} \,\frac{\partial \gamma }{\partial s}-\frac{u_{\tau}(\tau_{1}-\tau_{2})}{u}
\notag\\
 &=(\tau_{1}-\tau_{2})
  \left(
   \frac{\nabla u}{\sqrt{2}u}
   \cdot
   \frac{\sqrt{2}\,\gamma^{\prime}(s)}{\tau_{1}-\tau_{2}}-\frac{u_{\tau}}{u}
  \right)
\notag\\
& \leq (\tau_{1}-\tau_{2})
  \left(
   \frac{2|\gamma^{\prime}(s)|^2}{(\tau_{1}-\tau_{2})^2}
    + \frac{|\nabla u|^2}{2u^2}-\frac{u_{\tau}}{u}
  \right)
\notag\\
&= \frac{2|\gamma^{\prime}(s)|^2}{(\tau_{1}-\tau_{2})}
 +\frac{\tau_{1}-\tau_{2}}{2}
  \left(
       \frac{|\nabla u|^2}{u^2}-\frac{2u_{\tau}}{u} \,.
  \right)
\end{align}
By our gradient estimate, if $R\geq 0$, then
\begin{equation}
 \frac{|\nabla u|^2}{u^2}
 -\frac{2u_{\tau}}{u}\leq R+\frac{2n}{\tau}
\notag
\end{equation}
where $\tau=\tau_{2}+(1-s)(\tau_{1}-\tau_{2})$. Therefore
\begin{align}
 \frac{\partial \ell(s)}{\partial s}
 \leq
 \frac{2|\gamma^{\prime}(s)|^2}{(\tau_{1}-\tau_{2})}
  +\frac{\tau_{1}-\tau_{2}}{2}
  \left(
       R+\frac{2n}{\tau}
  \right) \,.
\end{align}
Integrating with respect to $s$ on $[0,1]$, we have
\begin{align}
 \ell(1)-\ell(0)
 \leq
 \frac{2\int_{0}^{1}|\gamma^{\prime}(s)|^2\,ds}{(\tau_{1}-\tau_{2})}
  +\frac{(\tau_{1}-\tau_{2})\int_{0}^{1}R\,ds}{2}
  +n\ln \frac{\tau_{1}}{\tau_{2}} \,.
\end{align}
Recall $\ell(0)=\ln\,u(x, t_1)$, $\ell(1)=\ln\,u(y, t_2)$, then
\begin{align}
\ln\frac{u(y, t_2)}{u(x, t_1)}
 \leq
 \frac{\int_{0}^{1} [\,4|\gamma^{\prime}(s)|^2+(\tau_{1}-\tau_{2})^2\,R\,]\,ds}{2(\tau_{1}-\tau_{2})}
  +\ln \left(
            \frac{\tau_{1}}{\tau_{2}}
       \right)^n \,.
\end{align}
Therefore, given any two points $(x, t_1)$, $(y, t_2)$ in the
space-time, we have
\begin{align}
u(y, t_2)
 \leq
 u(x, t_1)
 \left(
            \frac{\tau_{1}}{\tau_{2}}
       \right)^n
 \exp
  {
   \frac
        {\int_{0}^{1}      [\,4|\gamma^{\prime}(s)|^2+(\tau_{1}-\tau_{2})^2\,R\,]\,    ds}
        {2(\tau_{1}-\tau_{2})} \,.
  }
\end{align}

\end{proof}


\section{Localized version of the Gradient Estimate in section 2}


In this section we prove a localized version of the previous
gradient estimate. Here we apply Li-Yau's idea of using certain
cut-off functions to the new equations derived in the last
section. However the computation is more complicated for two
reasons. One is that  the metric is also evolving. The other is
that the equations coming from the last section have a more
complex structure.

\begin{theorem}Let M be a compact Riemannian manifold equipped with a family of
Riemannian metrics evolving under Ricci flow, that is,
$\frac{\partial g_{ij}}{\partial t}=-2R_{ij}$. Given $x_0 \in {\bf
M}$ and $r_0>0$, let $u$ be a smooth positive solution to the
conjugate heat equation $\Box^{*}u=-\triangle u-u_{t}+Ru=0$ in the
cube $Q_{r_0,T}:=
         \{(x,t)\,| d(x,x_{0},t) \leq r_0, 0<t \le T \}$ and
          $\tau=T-t$. Suppose $Ric\geq -K$ throughout $Q_{r_0,T}$ for some positive constant $K$. Then
for $(x,t)$ in the half cube $Q_{\frac{r_{0}}{2},\frac{T}{2}}:=
         \{(x,t)\,| d(x,x_{0},t)\leq \frac{r_{0}}{2}, 0<t\leq
         \frac{T}{2}\}$,
we have
\begin{align}
\frac{|\nabla u|^2}{u^2}-\frac{2u_{\tau}}{u}-R \leq
   c\,K
   +\frac{c}{T}
   +\frac{c}{r_{0}^2}
\end{align}
where $c>0$ is a constant depending only on the dimension $n$.
\end{theorem}
\medskip

\begin{proof}

As before let $f$ be a function defined by
$u=\frac{e^{-f}}{(4\pi\tau)^\frac{n}{2}}$ and
\[
q=2\triangle f-|\nabla f|^2+R = \frac{|\nabla
u|^2}{u^2}-\frac{2u_{\tau}}{u}-R.
\]
 From
inequality (\ref{heat-operator-q}) in the last section, we have
\begin{align}
    \triangle  q - q_{\tau}
    -2\nabla f \,\nabla q
    \geq \frac{1}{2n}\big(q+|\nabla f|^2+R\big)^2 \,.
\end{align}
For the fixed point $x_{0}$ in $M$, let $\varphi(x,t)$ be a smooth
cut-off function(mollifier) with support in the cube
\begin{align}
Q_{r_{0},T}:=\{(x,t)\,|\,x\in M, d(x,x_{0},t)\leq r_{0}, 0<t\leq
T\}
\end{align}
possessing the following properties:
\begin{align}
 & (1)\quad \varphi=\varphi(d(x,x_{0},t),t) \equiv
 \psi(r(x,t))\,\eta(t), r(x,t)=d(x,x_{0},t);\, \frac{\partial \psi}{\partial r}\leq0,\frac{\partial \eta}{\partial t}\leq0, \tau=T-t;
\notag\\
 & (2)\quad \varphi(x,t)\equiv1\,\,in\,\,
     Q_{\frac{r_{0}}{2},\frac{T}{2}}:=
         \{(x,t)\,| d(x,x_{0},t)\leq \frac{r_{0}}{2}, 0<t\leq \frac{T}{2}\};
\notag\\
 & (3)\quad |\frac{\partial_{r} \psi}{{\psi}^a}| \leq
 \frac{c(n,a)}{r_{0}},\,
   |\frac{\partial_{rr} \psi}{{\psi}^a}| \leq
   \frac{c(n,a)}{r_{0}^2},\, \text{for some $c(n,a)$, $0 < a< 1$};
\notag\\
& (4)\quad |\frac{\partial_{t} \eta}{\sqrt{\eta}}| \leq
 \frac{c}{T},
 \text{ for some $c$ depending on n}\notag \,.
\end{align}
Now we focus on the product $(\varphi\,q)(x,t)$. Since $\varphi$ has
support in $Q_{r_{0}, T}$, we can assume $\varphi\, q$ reaches its
maximum at some point $(y,s)\in Q_{r_{0},T}$. If $q(y,s)=2\triangle
f(y,s)-|\nabla f(y,s)|^2+R(y,s)$ is negative, then the theorem is
trivially true. Thus we can assume $q(y,s)\geq 0$. By direct
computation,
\begin{align}  \label{heat-operator-varphi-q}
\quad\triangle  (\varphi \,q) - (\varphi\,q)_{\tau}
    -2\nabla f \,\nabla (\varphi \,q)-2\frac{\nabla
    \varphi}{\varphi}\,\nabla(\varphi\,q)
\end{align}
\begin{align}
& =\varphi \big(\triangle q-q_{\tau}-2\nabla f\nabla q\big)
    +(\triangle \varphi) q -2\frac{|\nabla \varphi|^2}{\varphi}\,q
    -q\,\varphi_{\tau}-2q\nabla f\nabla \varphi
\notag\\
    &\geq \frac{\varphi}{2n}\big(q+|\nabla f|^2+R\big)^2
    +(\triangle \varphi) q -2\frac{|\nabla \varphi|^2}{\varphi}\,q
    -q\,\varphi_{\tau}-2q\nabla f\nabla \varphi.
\notag
\end{align}
At the point $(y,s)$ where the maximum value for $\varphi q$ is
attained, there hold
\begin{align}
 &(1)\quad\triangle(\varphi q)(y,s)\leq 0\notag\\
 &(2)\quad\nabla(\varphi q)(y,s)= 0\notag\\
 &(3)\quad(\varphi q)_{\tau}(y,s)\geq 0\notag \,.
\end{align}
The last inequality comes from the fact that
$(\varphi\,q)|_{t=T}=0$ since $\varphi|_{t=T}=0$, $\varphi\,q$ can
only take its maximum for $t \in [0,T)$. We have also borrowed the
idea of Calabi as used in \cite{LY} to circumvent the possibility
that $(y, s)$ is in the cut locus of $g(s)$.

Thus at the point $(y,s)$, inequality (\ref{heat-operator-varphi-q})
becomes
\begin{align}
    \frac{\varphi}{2n}\big(q+|\nabla f|^2+R\big)^2(y,s)
    &\leq
    \,\varphi_{\tau}q+2q\nabla f\nabla \varphi+2\frac{|\nabla \varphi|^2}{\varphi}\,q
    -(\triangle \varphi) q
\\
&=(I)+(II)+(III)+(IV) \notag  \,.
\end{align}
We estimate each term on the right-hand side by the following:
\begin{flushleft}
$(I)=\psi_{\tau}\,\eta(\tau)q+\psi \,\eta_{\tau}\, q
    = \psi_{r}\,r_{\tau}\,\eta(\tau)\,q+\psi \,\eta_{\tau}\,q$\,.
\end{flushleft}
From the lower bound assumption on the Ricci curvature
$Ric\geq-K$, we have (see [CK] e.g.)
\begin{align}
\frac{\partial r}{\partial \tau}=-\frac{\partial r}{\partial
t}\geq-K\,r  \,.
\end{align}
By construction of $\psi$, we have $\psi_{r}\leq0$. Therefore,
\begin{align}
\psi_{r}\,\frac{\partial r}{\partial \tau}
 \leq
 K\,r|\psi_{r}|= K\,r
 \frac{|\psi_{r}|}{\sqrt{\psi}}\,\sqrt{\psi}
 \leq
 K\,r\,\frac{c}{r_{0}}\sqrt{\psi}
 \leq
 c\,K\sqrt{\psi}  \,.
\end{align}
This shows
\begin{flushleft}
$(I)\leq c\,K \sqrt{\psi}\,\eta\, q+\psi \eta_{\tau}q
 \leq c\,K \sqrt{\varphi}\, q+ \psi \frac{|\eta_{\tau}|}{\sqrt{\eta}}\sqrt{\eta}\,q
 \leq c\,K \sqrt{\varphi}\, q+ \sqrt{\psi}\frac{c}{T}\sqrt{\eta}q. $
\end{flushleft}

Recall that $\varphi = \psi \eta$,  for a parameter $\epsilon$ to
be chosen later, we have
\begin{align}
(I)\leq c\,K \sqrt{\varphi}\, q+ \frac{c}{T}\sqrt{\varphi}\,q
 \leq
 \frac{c}{\epsilon}\,K^2+\frac{c}{T^2}+2\,\epsilon\,\varphi\,q^2 \,.
\end{align}
%
\begin{align}
(II)
&\leq
 2|\nabla f |\,|\nabla\, \varphi|\,q
 \leq
 2|\nabla f |\,\frac{|\nabla\, \varphi|}{\sqrt{\varphi}}\,\sqrt{\varphi}\,q
\notag\\
&\leq
 \frac{1}{\epsilon}\,4\,|\nabla f |^2 \frac{|\nabla \varphi|^2}{\varphi}
 +\epsilon \,\varphi\, q^2
\notag\\
\qquad
&=
 \frac{1}{\epsilon}\,4\,\sqrt{\varphi}\,|\nabla f |^2\, \frac{|\nabla
 \varphi|^2}{\varphi^{\frac{3}{2}}}
 +\epsilon \,\varphi\, q^2
\notag\\
 &\leq
\epsilon\, \varphi\, |\nabla f |^4
  +\frac{1}{\epsilon^{3}} \frac{c}{r_{0}^2}
  + \epsilon\,\varphi\, q^2  \,.
\end{align}
\begin{align}
(III)
 =2\,\frac{|\nabla \varphi|^2}{\varphi^{3/2}}\,\sqrt{\varphi}\,q
 \leq 2\,|\frac{\nabla \varphi}{\varphi^{3/4}}|^2\,\sqrt{\varphi}\,q
 \leq \epsilon \varphi \,q^2+\frac{1}{\epsilon}\frac{c}{r_{0}^2} \,.
\end{align}
\vspace{-12pt}
\begin{align}
&\negthickspace\negthickspace\negthickspace\negthickspace
\negthickspace\negthickspace\negthickspace\negthickspace
\negthickspace\negthickspace\negthickspace\negthickspace
\negthickspace\negthickspace\negthickspace\negthickspace
\negthickspace\negthickspace
(IV)
 =-(\triangle \varphi)\,q
 =
     -\left(
       \partial_{rr}\varphi
       +(n-1)\frac{\partial_{r}\varphi}{r}
       +\partial_{r}\varphi\,\partial_{r}\log \sqrt{g}
       \right)\,q
\notag\\
 & \qquad
  \leq
  \frac{|\partial_{rr}\varphi|}{\sqrt{\varphi}}\, \sqrt{\varphi}\,q
  +(n-1)\frac{\partial_{r}\varphi}{r}\,q
  +\frac{|\partial_{r}\varphi|}{\sqrt{\mathstrut \varphi}} \sqrt{\mathstrut K}\,\sqrt{\mathstrut \varphi} \,q
\notag\\
 & \qquad
  \leq
      \frac{1}{\epsilon}
      \left(
            \frac{\partial_{rr}\varphi}{\sqrt{\varphi}}
      \right)^2
      +\epsilon \,\varphi\,q^2
      +(n-1)\frac{\partial_{r}\varphi}{r}\,q
      +\frac{K}{\epsilon}
      \left(
             \frac{|\partial_{r} \varphi|}{\sqrt{\varphi}}
       \right)^2
      +\epsilon \,\varphi\,q^2
\notag\\
 & \qquad
  \leq
      2\epsilon \,\varphi\,q^2
      +\frac{c}{r_{0}^4}+\frac{c\,K}{r_{0}^2}
      +(n-1)\frac{\partial_{r}\varphi}{r}\,q  \,.
\notag
\end{align}
Notice $\varphi\equiv 1$ in $Q_{\frac{r_{0}}{2},\frac{T}{2}}$,
thus $\partial_{r}\varphi=0 $ for $0\leq r\leq \frac{r_{0}}{2}$,
we can just focus on $r\geq \frac{r_{0}}{2}$, $\frac{1}{r}\leq
\frac{2}{r_{0}}$, then $(IV)$ can be estimated as
\begin{align}
&(IV)
  \leq
      2\epsilon \,\varphi\,q^2
      +\frac{c}{r_{0}^4}+\frac{c\,K}{r_{0}^2}
      +(n-1)\frac{2\partial_{r}\varphi}{r_{0}}\,q
\notag\\
 & \qquad
  \leq
      2\epsilon \,\varphi\,q^2
      +\frac{c}{r_{0}^4}+\frac{c\,K}{r_{0}^2}
      +\epsilon
      \,\varphi\,q^2+\frac{4(n-1)^2}{r_{0}^2}|\frac{\partial_{r}\varphi}{\sqrt{\varphi}}|^2
\notag\\
 & \qquad
  \leq
      3\epsilon \,\varphi\,q^2
      +\frac{c}{r_{0}^4}+\frac{c\,K}{r_{0}^2} \,.
\notag
\end{align}
Combine (I)-(IV), we have
\begin{align}  \label{H>=2nd}
 \frac{\varphi}{2n}\big(q+|\nabla f|^2+R\big)^2(y,s)
 & \leq
  7\epsilon \varphi\,q^2
  +c\,K^2
  +\frac{c}{T^2}
  +\frac{c}{r_{o}^2}
  +\frac{cK}{r_{o}^2}
  +\frac{c}{r_{o}^4}
  +\epsilon\, \varphi\,|\nabla f|^4
\notag\\
 & \leq
  7\epsilon \varphi\,q^2
  +c\,K^2
  +\frac{c}{T^2}
  +\frac{c}{r_{o}^4}
  +\epsilon\, \varphi\,|\nabla f|^4 \,.
\end{align}
Notice we assumed $q(y,s)\geq 0$, otherwise the theorem is trivially
true,
\begin{align}
 \big(q+|\nabla f|^2+R\big)^2(y,s)
 &=\big(q+|\nabla f|^2+R^{+}-R^{-}\big)^2(y,s)
\notag\\
 & \geq
  \frac{1}{2}\big(q+|\nabla f|^2+R^{+}\big)^2(y,s)-(R^{-})^{2}(y,s)
\notag\\
& \geq
  \frac{1}{2}\big(q+|\nabla f|^2\big)^2(y,s)-(R^{-})^{2}(y,s)
\notag\\
& \geq
  \frac{1}{2}\big(q^2+|\nabla f|^4\big)(y,s)-(\sup_{Q_{r_{0},T}} R^{-})^{2}
\notag\\
& \geq
  \frac{1}{2}\big(q^2+|\nabla f|^4\big)(y,s)-n^2K^2.
\end{align}
Here we have used the inequalities $2(a-b)^2\geq a^2-2b^2$,
$(a+b)^2\geq a^2+b^2$ for $a,\,b\geq 0$ and the lower bound
assumption for the Ricci curvature $Ric\geq -K\,\,\Rightarrow
\,\,R\geq -nK\,\,\Rightarrow \,\,R^{-}\leq nK$ since $R=-R^{-}$ if
$R<0$. Substituting into (\ref{H>=2nd}) and reorganizing, we have
\begin{align}
 (\frac{1}{4n}-7\epsilon)\varphi \,q^2(y,s)
  \leq
  (\epsilon -\frac{1}{4n})\varphi |\nabla f|^4
  +c\,K^2
  +\frac{c}{ T^2}
  +\frac{c}{ r_{o}^4}.
\end{align}
Take $\epsilon$ such that $7\epsilon \leq\frac{1}{4n}$, then the
above inequality becomes
\begin{align}
  \varphi \,q^2(y,s)
  \leq
  c\,K^2
  +\frac{c}{T^2}
  +\frac{c}{r_{o}^4} \,.
\end{align}
By using inequality $a_{1}^2+a_{2}^2+...+a_{n}^2\leq
(a_{1}+a_{2}+...+a_{n})^2$,
\begin{align}
 (\varphi \,q)^2(y,s)
 \leq
  \varphi \,q^2(y,s)
 \leq
  \big(
   c\,K
   +\frac{c}{T}
   +\frac{c}{r_{o}^2}
  \big)^2 \,.
\end{align}
Since if $(x,t)\in Q_{\frac{r_{o}}{2},\frac{T}{2}}$, then
$\varphi(x,t)\equiv1$, thus for any $(x,t)\in
Q_{\frac{r_{o}}{2},\frac{T}{2}}$,
\begin{align}
q(x,t)&=\varphi(x,t)\,q(x,t)\leq
\max_{Q_{\frac{r_{o}}{2},\frac{T}{2}}}(\varphi \,q)(x,t) \leq
\max_{Q_{r_{o},T}}(\varphi \,q)(x,t)=(\varphi \,q)(y,s) \notag\\
& \leq c\,K
   +\frac{c}{T}
   +\frac{c}{r_{o}^2} \,.
\end{align}
Therefore we just proved that in
$Q_{\frac{r_{o}}{2},\frac{T}{2}}$,
\begin{align}
q(x,t) \leq
   c\,K
   +\frac{c}{T}
   +\frac{c}{r_{0}^2} \,.
\end{align}
If we bring back $u$, recall $q=2\triangle f-|\nabla f|^2+R$,
$f=-\ln u-\frac{n}{2}\ln (4\pi\tau)$, then we have
\begin{align}
\frac{|\nabla u|^2}{u^2}-\frac{2u_{\tau}}{u}-R \leq
   c\,K
   +\frac{c}{T}
   +\frac{c}{r_{o}^2}  \,.
\end{align}
\end{proof}

\section{A Generalization of Perelman's W-Entropy}

In this section, we show that Perelman's W-entropy and its
monotonicity can be generalized to a wider class. It is an
established fact that monotonicity formulas tend to provide useful
information on the underlining equation. Therefore the more of them
are found the better. For a related but different generalization of
Perelman's formula and its applications, please see a recent paper
\cite{Li}.

Define a family of  entropy formulas for the Ricci flow case by:
\begin{equation}
\label{newentropy}
W(g,f,\tau):=\int_{M}
                   \left(
                   \frac{a^2}{2\pi}\tau(R+|\nabla f|^2)+f-n
                   \right)u\,dx
\end{equation}
where $R$ is the scalar curvature, $\tau=T-t>0$;
$u=\frac{e^{-f}}{(4\pi\tau)^\frac{n}{2}}$ is a positive solution to
the following conjugate heat equation (\ref{conjugate}), satisfying
$\int u \,dx=1$,
\begin{equation}  \label{conjugate}
\Box^*u=-u_{t}-\triangle u +R\,u=0 \,.
\end{equation}
Notice that
\begin{equation}
  \begin{cases}
    f=-\ln u-\frac{n}{2}\ln(4\pi\tau)  \\
    \nabla f=-\frac{\nabla u }{u} \\
    \triangle f=-\frac{\triangle u }{u}+|\nabla f|^2  \\
    \frac{\partial f}{\partial t}=-\frac{u_{t}}{u}+\frac{n}{2\tau}
  \end{cases}
\end{equation}
we get the evolution equation for $f$,
\begin{equation} \label{evolution-f}
   \frac{\partial f}{\partial t}
    =-\triangle f +|\nabla f|^2-R+\frac{n}{2\tau}  \,.
\end{equation}
Now we come to the theorem of the section:

\begin{theorem}
\label{Ricci} Let g(t) evolve by the Ricci flow, that is,
$\frac{\partial g}{\partial t}=-2 Ric $ on a closed manifold $M$ for
$t\in [0,T)$, and $u:M\times[0,T)\mapsto(0,\infty)$ with
$u=\frac{e^{-f}}{(4\pi\tau)^\frac{n}{2}}$ be a positive solution to
the conjugate heat equation (\ref{conjugate}). For $0\leq
a^2\leq2\pi$,  the functional defined in (\ref{newentropy}) is
increasing according to
\begin{equation}
 \frac{\partial }{\partial t} W(g,f,\tau)
 \geq\frac{a^2\tau}{\pi} \int_{M}
 |Ric+Hess(f)-\frac{g}{2\tau}|^2u\,dx\geq 0.
\end{equation}
\end{theorem}

We start the proof with the following results by Perelman \cite{P}
in the form of lemmas from [T]. For completeness, the proofs are
given in the Appendix.
\begin{lemma}
\label{Ricci-F} Let $u,f,g,\tau$ defined as in Theorem
\ref{Ricci}, then the F-entropy defined by
 $F(g,f):=\int \big(R+|\nabla f|^2\big)u\,dx$ is non-decreasing in
 $t$ under
\begin{equation}
\frac{\partial F(g,f)}{\partial t}=2\int |Ric+Hess(f)|^2 u\,dx
\end{equation}
\end{lemma}
\begin{lemma}  \label{Ricci-P}
Let $g,f,\tau$ defined as in Theorem \ref{Ricci}, define
 $P:=[\tau(2\triangle f-|\nabla f|^2+R)+f-n]u$, then
\begin{equation}
\Box^{*}P=-2\tau|Ric+Hess(f)-\frac{g}{2\tau}|^2u  \,.
\end{equation}
\end{lemma}

\begin{proof} (of Theorem \ref{Ricci})
Notice
\begin{align}
  W(g,f,\tau)
  &=\,\int_{M}
                   \left(
                   \frac{a^2}{2\pi}\tau(R+|\nabla f|^2)+f-n
                   \right)u\,dx
  \notag\\
  &=\frac{a^2}{2\pi}\int_{M}
                   [
                   \tau(R+|\nabla f|^2)+f-n
                   ]u\,dx
 \notag
 \\
 &\quad+\big(1-\frac{a^2}{2\pi}\big)
      \left(
      \int_{M}f\,u\,dx
      \right)-\big(1-\frac{a^2}{2\pi}\big)n
\notag
\end{align}
we split the derivative of $W$ over the time $t$ into two parts,
\begin{align}
 \frac{\partial }{\partial t} W(g,f,\tau)
 &= \frac{a^2}{2\pi}\frac{\partial }{\partial t}
            \left(   \int_{M}
                      [
                      \tau(R+|\nabla f|^2)+f-n
                      ]
                      u\,dx
            \right)
 \notag
 \\
 &\quad+\big(1-\frac{a^2}{2\pi}\big)\frac{\partial }{\partial t}
      \left(
      \int_{M}f\,u\,dx
      \right)
\notag\\
&=\frac{a^2}{2\pi}\frac{\partial }{\partial t}\int_{M}P\,dx
   +\big(1-\frac{a^2}{2\pi}\big)\frac{\partial }{\partial t}
      \left(
      \int_{M}f\,u\,dx
      \right) \,.
\notag
\end{align}
We compute for each term,
\begin{align}
 \frac{a^2}{2\pi}\frac{\partial }{\partial t}\int_{M}P\,dx
 &=\frac{a^2}{2\pi} \int_{M}
   \big(
     P_{t}\,dx
   + P\frac{\partial dx}{\partial t}
   \big)
 =\frac{a^2}{2\pi} \int_{M}
   \big(
    P_{t}\,dx
   + P(-R)\,dx
   \big)
 \notag\\
 &=-\frac{a^2}{2\pi} \int_{M} \Box^{*}P\,dx
   -\frac{a^2}{2\pi} \int_{M} \triangle P\,dx
 =-\frac{a^2}{2\pi} \int_{M} \Box^{*}P\,dx
 \notag
\end{align}
the last equality comes from $\int_{M}\triangle P\,dx =0$ for closed
manifold $M$. By Lemma \ref{Ricci-P}, we have
\begin{align}
 \frac{a^2}{2\pi}\frac{\partial }{\partial t}\int_{M}P\,dx
 &
 =
 \frac{a^2\tau}{\pi}
 \int_{M} |Ric+Hess(f)-\frac{g}{2\tau}|^2u\,dx\geq 0 \,.
 \notag
\end{align}
It suffices to prove the non-negativity of $\frac{\partial
}{\partial t}\left(\int_{M}f\,u\,dx\right)$. Follow the direct
computation,
\begin{align}
 \frac{\partial }{\partial t}\left(\int_{M}f\,u\,dx\right)
 &=\int_{M} f_{t}\,u\,dx +f\,u_{t}\,dx+f\,u\, \frac{\partial (dx)}{\partial t}
\notag\\
 &=\int_{M} (-\triangle f+|\nabla f|^2-R+\frac{n}{2\tau})u\,dx
 \notag\\
 &\quad
 +\int_{M}
   \big(
   f(-\triangle u+Ru)-Rfu
   \big)
   \,dx  \,.
 \notag
\end{align}
Using integration by parts, we have
\begin{align}  \label{evolve-F}
 \frac{\partial }{\partial t}\left(\int_{M}f\,u\,dx\right)
 &=\int_{M} (-2\triangle f+|\nabla f|^2)u\,dx +\int_{M} (\frac{n}{2\tau}-R)u\,dx
\notag\\
 &=\int_{M} (2\frac{\triangle u}{u} -|\nabla f|^2)u\,dx +\int_{M} (\frac{n}{2\tau}-R)u\,dx
\notag\\
 &=\int_{M} -|\nabla f|^2\,u\,dx +\int_{M} (\frac{n}{2\tau}-R)u\,dx
\notag\\
 &=\frac{n}{2\tau}-\int_{M} (|\nabla f|^2+R)\,u\,dx \,.
\end{align}
Now we turn to estimate of $F(g,\tau)=\int_{M}(|\nabla
f|^2+R)u\,dx$.\\ From Lemma \ref{Ricci-F}, we have
\begin{align}
\frac{\partial F}{\partial t}
 &=2\int |Ric+Hess(f)|^2u\,dx=2\int \big(\sum_{i,j} |R_{ij}+f_{ij}|^2\big)u\,dx
\notag\\
 &\geq 2\int \big(\sum_{i=j} |R_{ij}+f_{ij}|^2\big)u\,dx
  \geq 2\int \frac{1}{n}\big(\sum R_{ii}+\sum f_{ii}\big)^2 u\,dx
\notag\\
 &= \frac{2}{n}\int \left(R+\triangle f \right)^2 u\,dx  \,.
\end{align}
The last inequality comes from
$\sqrt{\frac{a_{1}^2+...+a_{n}^2}{n}}\geq
\frac{a_{1}+...+a_{n}}{n}$ for $a_{i}\geq 0$. Also by
Cauchy-Schwarz inequality, we have
\begin{align}
 \int (R+ \triangle f) \sqrt{u}\,\sqrt{u}\,dx
\leq
 \left(\int (R+\triangle f)^2 u\,dx\right)^{\frac{1}{2}}
 \left(\int u\,dx\right)^{\frac{1}{2}}  \,.
\end{align}
Since $\int u\,dx=1$, the above inequality can be simplified as
\begin{align}
 \left(\int (R+\triangle f) u\,dx\right)^2
\leq \int (R+\triangle f)^2 u\,dx\,.
\end{align}
Then the evolution of $F$ along the time $t$ would be estimated by
\begin{align}
\frac{\partial F}{\partial t}
 &\geq
\frac{2}{n} \left(\int (R+\triangle f) u\,dx\right)^2
 =\frac{2}{n} \left(\int (R+|\nabla f|^2) u\,dx\right)^2
\end{align}
due to the following equality in closed manifold $M$
\begin{align}
\int_{M} (\triangle f-|\nabla f|^2) u\,dx
 &=\int_{M} (-\frac{\triangle u}{u}+|\nabla f|^2-|\nabla f|^2) u\,dx
 \notag\\
 & =-\int_{M} \triangle u\,dx=0
 \notag
\end{align}
\begin{align}
\Rightarrow
 \,\,\int_{M} (\triangle f) \,u\,dx= \int_{M}|\nabla f|^2 u\,dx  \,.
\end{align}
From the definition $F=\int (R+|\nabla f|^2) u\,dx$, we get
\begin{align}
\frac{\partial F}{\partial t}
 &\geq
 \frac{2}{n}F^2\geq 0  \,.
\end{align}
We claim
\begin{align}
F(t)\leq
 \frac{n}{2(T-t)}  \,.
\end{align}
Here is the proof of the above claim,
\begin{align}
&\frac{dF}{dt}
 \geq
 \frac{2}{n}F^2
\,\,\Rightarrow \,\,
 \frac{dF}{F^2}
 \geq
 \frac{2}{n}dt
\,\,\Rightarrow \,\,
 \int_{t}^{T}
\frac{dF}{F^2}
 \geq
 \frac{2}{n}(T-t)
 \notag\\
\Rightarrow \,\,&\,\,
 -(\frac{1}{F(T)}-\frac{1}{F(t)})
 \geq
\frac{2}{n}(T-t) \,\,\Rightarrow \,\, \frac{1}{F(t)}
 \geq
\frac{2}{n}(T-t)+\frac{1}{F(T)} \notag  \,.
\end{align}
If $F(T)> 0$, then
 $\frac{1}{F(t)} \geq \frac{2}{n}(T-t)$, that is,
 $F(t)\leq \frac{n}{2(T-t)}$;\\
If $F(T)\leq 0$, since $\frac{dF}{dt}\geq 0$, then $F(t)\leq 0\leq
\frac{n}{2(T-t)}$ for all $t\in[0,T)$, therefore,
\begin{align}
F(t)=\int (R+|\nabla f|^2) u\,dx
 \leq
 \frac{n}{2(T-t)}=\frac{n}{2\tau}
\end{align}
plugging into (\ref{evolve-F}), we obtain
\begin{align}
 \frac{\partial }{\partial t}\left(\int_{M}f\,u\,dx\right)
=\frac{n}{2\tau}-\int_{M} (|\nabla f|^2+R)\,u\,dx\geq 0  \,.
 \notag
\end{align}
Thus we complete the proof of Theorem \ref{Ricci}.

\end{proof}

\section{The case for the heat equation under a fixed metric}

It is well known that gradient estimates and monotonicity formulas
for the heat equation in the fixed metric case are important in
their own right. The Li-Yau estimate is one of several examples.
Recently Perelman's W-entropy and gradient estimate for the Ricci
flow were transformed to the case of the heat equation in the
fixed metric case by Lei Ni \cite{N1}. He also pointed out some
useful geometric applications.

Here we will transplant some of the results in the previous
sections to the heat equation case. More specifically, we will
introduce a family of entropy formulas which include both
Perelman's
 $W$-entropy (for the heat equation as defined in \cite{N1}) and the 'Boltzmann-Shannon' entropy
 as considered in thermodynamics, information theory.

\indent Let $(M,g)$ be a closed Riemannian $n$-manifold with the
metric $g$ not evolving along time $t$, suppose $u$ is a positive
solution to the heat equation (\ref{heat}) with $\int_{M}udx=1$,
\begin{equation}
(\triangle-\frac{\partial}{\partial t})u=0 \label{heat}  \,.
\end{equation}
Let $f:M\rightarrow R$ be smooth, defined by
  $f:=-\ln u-\frac{n}{2}\ln(4\pi \tau)$,
that is,
  $u=\frac{e^{-f}}{(4\pi\tau)^{\frac{n}{2}}}$, where $\tau=\tau(t)>0$ is a scale parameter
with $\frac{d\tau}{dt}=1$.  We define a family of entropy formulas
by:
\begin{equation}
W(f,\tau,a):=\int_{M}
 \left(
       \frac {a^2 \tau}
             {2\pi}
       |\nabla f|^2 +f - n
      +\frac{n}{2}
                 \ln \frac {2\pi}
                           {a^2}
 \right)
 \frac {e^{-f}}
       {(4\pi \tau)^\frac{n}{2} }
 dx  \label{entropy}  \,.
\end{equation}

The main results of the section are the next two theorems.

\begin{theorem}  \label{maintheorem-1}
\indent Let M be a closed Riemannian $n$-manifold with fixed metric
$g$, and $u$, $f$ be defined as above satisfying $\int_{M}udx=1$ and
$(\triangle-\frac{\partial}{\partial t})u=0$. The entropy
$W(f,\tau,a)$ defined in (\ref{entropy}) satisfies:
\begin{align}
 \frac {\partial W(f,\tau,a)}
       {\partial t}
   & =-\frac{a^2 \tau}{\pi}
       \int_{M} u
               |f_{ij}-\frac{g_{ij}}{2\tau}|^2dx
        -\frac{a^2\tau}{\pi}\int_{M}u R_{ij}f_{i}f_{j}\,dx
        \\
   &\,\,\,\,\,\,\,
    + \big(1-\frac{a^2}{2\pi}\big)
                     \int_{M} u
                               \left(
                                     \frac {|\nabla u|^2}
                                           {u^2}
                                     -\frac{\triangle u}{u}
                                     -\frac{n}{2\tau}
                               \right)dx \notag  \,.
\label{monotone}
\end{align}

Moreover, if $0\leq a^2\leq 2\pi$, $\tau=t$ and $M$ has non-negative
Ricci curvature, then $W$ is monotone non-increasing in time t, that
is,
\begin{equation}
 \frac {\partial W(f,t,a)}
       {\partial t}\leq 0 \,.
\end{equation}
\end{theorem}
\begin{remark}
If we choose $a^2=2\pi$ in $W(f,\tau,a)$, we recover Perelman's
$W$-entropy worked out in \cite{N1} and its monotonicity as a
special case.
\begin{equation}
W(f,\tau,\sqrt{2\pi})=\int_{M}
  \left(
  \tau |\nabla f|^2+f-n
  \right)
  u dx  \label{leini}
\end{equation}
\begin{equation}
\frac {\partial W(f,\tau,\sqrt{2\pi})}
      {\partial t}
      = - 2\tau \int_{M} u
            \left(
                  |f_{ij}- \frac
                  {g_{ij}}{2\tau}|^2+R_{ij}f_{i} f_{j}
            \right)
                   dx
      \leq 0 \label{perelman-monotone}  \,.
\end{equation}

\end{remark}

\begin{remark}
If $a=0$, $W(f,\tau,0)$ is related with 'Boltzmann-Shannon'
entropy as considered in thermodynamics, information theory,
\begin{equation}
 N = \int_{M} -u\, \ln u \,\,dx   \,.
 \label{shannon}
\end{equation}
The  entropy
 $W(f,\tau,a)$ thus defined serves as a connection between 'Boltzmann(1870s)-Shannon(1940s)' entropy and Perelman's
$W$-entropy(2002). We are partially motivated by the logarithm
Sobolev inequality of the Gross \cite{Gr} as stated in \cite{LL}.
\end{remark}

From this new entropy formula, we deduce the corresponding
differential inequality for the fundamental solution to the heat
equation.
\begin{theorem}  \label{maintheorem-2}
\indent Let M be a closed Riemannian $n$-manifold  with fixed metric
$g$ and non-negative Ricci curvature, $u$ be the fundamental
solution, and $f$ be defined as $f:=-\ln u-\frac{n}{2}\ln(4\pi t)$.
Then for any constant $\alpha \geq 1$,
\begin{equation}  \label{gradientestimate}
t\,
   \left(\alpha \,\triangle f \, -|\nabla f|^2
   \right)
+f-\alpha \,\frac{n}{2}\leq 0.
\end{equation}
\end{theorem}

\begin{remark}
In particular, if $\alpha=2$, then it becomes the following
differential inequality proved in [N1].
\begin{equation}
t\,(2 \triangle f \, -|\nabla f|^2)+f-n\leq 0  \,.
\end{equation}

If we divide the left-hand side of the inequality
(\ref{gradientestimate}) by $\alpha t$, where $\alpha \geq 1$, and
$t>0$, we get $ \triangle f \, -\frac{|\nabla
f|^2}{\alpha}+\frac{f}{\alpha\,t}-\frac{n}{2t}\leq 0$, let $\alpha
\rightarrow \infty$, then we conclude that the inequality
(\ref{gradientestimate}) includes Li-Yau gradient estimate, that
is, $\frac {|\nabla u|^2}{u^2}-\frac{\triangle
u}{u}-\frac{n}{2t}\leq 0$ since $\triangle f=\frac {|\nabla
u|^2}{u^2}-\frac{\triangle u}{u}$.
\end{remark}
\begin{remark}
For $\alpha > 2$, the gradient estimate (\ref{gradientestimate}) is
an interpolation of Perelman's gradient estimate cf. \cite{N1} and
Li-Yau estimate; \textbf{however}, for $1\leq \alpha \leq 2$, the
gradient estimate (\ref{gradientestimate}) is new here, it can't be
directly obtained from Perelman's gradient estimate and Li-Yau
gradient estimate.
\end{remark}
\begin{remark}
In the Euclidean Space $\mathbb{R}^n$, if $u$ is the fundamental
solution to the heat equation, then (\ref{gradientestimate})
becomes an equality.
\end{remark}

We start the proof of Theorem \ref{maintheorem-1} with the following
lemmas. Some computation results can be directly found in [LY], [N1]
or some other sources, we give the details here for completeness.

\begin{lemma}  \label{nabla-f}
Let $u$ be a positive solution to the heat equation (\ref{heat})
in a closed Riemannian n-manifold $M$, and $f=-\ln\,u
   -\frac{n}{2}
          \ln\,(4\pi\tau) $, where $\tau=\tau(t)> 0,$
           $\frac{\partial \tau}{\partial t}=1$. Then

\begin{equation}  \label{equa-nabla-f}
 (\triangle - \partial_{t}) (|\nabla f|^2)
 = 2 |f_{ij}|^2
   + 2\,\,\nabla f \,\,\nabla (|\nabla f|^2)
   +2 R_{ij}f_{i}f_{j}  \,.
\end{equation}
\end{lemma}

\begin{proof}
By direct computation,
\begin{align}
\triangle (|\nabla f|^2)
&=2 |f_{ij}|^2
   + 2( \,\,\partial_{j} f
   \,\,\partial_{i}\partial_{i} \partial_{j} f)
   \notag \\
&=2 |f_{ij}|^2
   + 2 \left(
             \,\,\partial_{j} f
   \,\,\partial_{j}\partial_{i}\partial_{i} f+ R_{ij}f_{i}f_{j}\,\,
   \right)
   \notag \\
&=2 |f_{ij}|^2
   + 2\,\nabla f \,\nabla (\triangle f)
   +2 \, R_{ij}f_{i}f_{j}\,.
\end{align}
Notice that
\begin{equation}
  \begin{cases}
    f=-\ln u-\frac{n}{2}\ln(4\pi\tau)  \\
    \triangle f=-\frac{\triangle u }{u}+|\nabla f|^2=-\frac{u_{t} }{u}+|\nabla f|^2  \\
    \frac{\partial f}{\partial t}=-\frac{u_{t}}{u}-\frac{n}{2\tau}
  \end{cases}
\end{equation}
then,
\begin{align}    \label{triangle-f}
\triangle f=f_{t}+|\nabla f|^2+\frac{n}{2\tau}\,.
\end{align}
Thus,
\begin{align}
\triangle (|\nabla f|^2) & =2 |f_{ij}|^2
   + 2\,\nabla f \,\nabla (f_{t}
   +|\nabla f|^2+\frac{n}{2\tau})
   +2 \, R_{ij}f_{i}f_{j}
   \notag\\
& =2 |f_{ij}|^2
   + 2 \,\nabla f \,\nabla f_{t}+2 \,\nabla f \,\nabla (|\nabla
   f|^2)
   +2 \, R_{ij}f_{i}f_{j}
   \notag\\
& =2 |f_{ij}|^2
   + (\nabla f|^2)_{t}+2 \,\nabla f \,\nabla (|\nabla
   f|^2)
   +2 \, R_{ij}f_{i}f_{j}
\end{align}
therefore,
\begin{equation}
(\triangle - \partial_{t}) (|\nabla f|^2) = 2 |f_{ij}|^2+
2\,\,\nabla f \,\,\nabla (|\nabla f|^2)+ 2 R_{ij}f_{i}f_{j}  \,.
 \notag
\end{equation}
\end{proof}

\begin{lemma} \label{triangle-f}
Let $u$,$f$ be defined as in Lemma \ref{nabla-f}, then
\begin{equation}   \label{equa-triangle-f}
(\triangle - \partial_{t}) (\,\triangle f) = 2 |f_{ij}|^2+
2\,\,\nabla f \,\,\nabla (\,\triangle f)+ 2 R_{ij}f_{i}f_{j}
\end{equation}
\begin{equation}
(\triangle - \partial_{t}) (2\,\triangle f-|\nabla f|^2) = 2
|f_{ij}|^2+ 2\,\,\nabla f \,\,\nabla (2\,\triangle f-|\nabla f|^2)+
2 R_{ij}f_{i}f_{j}   \,.
\end{equation}


\begin{proof}
By the result in Lemma \ref{nabla-f}, we have $ \triangle
f=f_{t}+|\nabla f|^2+\frac{n}{2\tau}$, then

\begin{align}
(\triangle - \partial_{t}) (\,\triangle f)
 & =\triangle(\,\triangle f)-\partial_{t}(\,\triangle f)
   =\triangle(f_{t}+|\nabla f|^2+\frac{n}{2\tau})-\partial_{t}(\,\triangle f)
\notag\\
 & =\triangle (f_{t})+\triangle (|\nabla f|^2)-\triangle (f_{t})
   =\triangle (|\nabla f|^2)
\notag\\
 & =2 |f_{ij}|^2
   + 2\,\,\nabla f \,\,\nabla (\,\triangle f)
   + 2 R_{ij}f_{i}f_{j}  \,.
\notag
\end{align}
Combining (\ref{equa-nabla-f}) and (\ref{equa-triangle-f}), we have,
\begin{equation}
(\triangle - \partial_{t}) (2\,\triangle f-|\nabla f|^2) =
2|f_{ij}|^2 + 2\,\,\nabla f \,\,\nabla(2\,\triangle f-|\nabla
f|^2)
 + 2 R_{ij}f_{i}f_{j}. \notag
\end{equation}

\end{proof}
\end{lemma}

\begin{lemma}    \label{F}
Let $u$,$f$ be defined as in Lemma \ref{nabla-f}, let $F=\frac
{a^2\tau} {2\pi}|\nabla f|^2 +f$, then
\begin{align}
(\triangle - \partial_{t}) F
&=\frac{a^2\tau}{\pi}|f_{ij}-\frac{g_{ij}}{2\tau}|^2
   +\big(1-\frac{a^2}{2\pi}\big)(|\nabla f|^2+ \frac{n}{2\tau})
\notag
\\
&\,\,\,\,\,\, +\frac{a^2}{\pi}\triangle f
   +\frac{a^2\tau}{\pi}\,\, \nabla f \,\,\,\nabla (|\nabla f|^2)
   +\frac{a^2\tau}{\pi}R_{ij}f_{i}f_{j}  \,.
\notag
\end{align}
\end{lemma}


\begin{proof}
Keep in mind that in local normal coordiantes $\sum
f_{ij}g_{ij}=\triangle f$, $\sum(g_{ij})^2=n$ and the result in
Lemma \ref{nabla-f}, we follow the direct computation,
\begin{align}
(\triangle - \partial_{t}) F
&=(\triangle-\partial_{t})(\frac{a^2\tau} {2\pi}|\nabla f|^2 +f)
   \notag \\
&=\frac{a^2\tau}{2\pi}(\triangle -\partial_{t})(|\nabla f|^2)
  -\frac{a^2}{2\pi}|\nabla f|^2+\triangle f-f_{t}
\notag\\
&=\frac{a^2\tau}{2\pi} \left(  2 (f_{ij})^2 + 2\,\,\nabla f
                               \,\,\nabla (|\nabla f|^2)+2R_{ij}f_{i}f_{j}
                       \right)
\notag\\
&\quad -\frac{a^2}{2\pi}|\nabla f|^2
   +|\nabla f|^2
   +\frac{n}{2\tau}
\notag\\
&=\frac{a^2\tau}{\pi}|f_{ij}-\frac{g_{ij}}{2\tau}|^2
   -\frac{a^2}{2\pi}\frac{n}{2\tau}
   +\frac{a^2}{\pi}\triangle f
   +\frac{a^2\tau}{\pi}\,\nabla f \,\nabla (|\nabla f|^2)
\notag\\
 &\quad
   +\frac{a^2\tau}{\pi}R_{ij}f_{i}f_{j}
   -\frac{a^2}{2\pi}|\nabla f|^2
   +|\nabla f|^2
   +\frac{n}{2\tau}
\notag
\\
&=\frac{a^2\tau}{\pi}|f_{ij}-\frac{g_{ij}}{2\tau}|^2
   +\big(1-\frac{a^2}{2\pi}\big)(|\nabla f|^2+ \frac{n}{2\tau})
\notag
\\
&\,\,\,\,\,\, +\frac{a^2}{\pi}\triangle f
   +\frac{a^2\tau}{\pi}\,\, \nabla f \,\,\,\nabla (|\nabla f|^2)
   +\frac{a^2\tau}{\pi}R_{ij}f_{i}f_{j}\,.
\notag
\end{align}
\end{proof}
\medskip

Now we turn to the
\begin{proof}
(of Theorem \ref{maintheorem-1})  We have,
\begin{align}
W(f,\tau,a)&=\int_{M}
 \left(
       \frac {a^2 \tau}
             {2\pi}
       |\nabla f|^2 +f
 \right)u dx- n
      +\frac{n}{2}
                 \ln \frac {2\pi}
                           {a^2}\notag\\
 &=\int_{M} Fu \,\, dx -n
      +\frac{n}{2}
                 \ln \frac {2\pi}
                           {a^2} \notag
\end{align}
then
\begin{align}
    \frac {\partial W}
          {\partial t}
    &= \int_{M} \frac{\partial}{\partial t}
            (Fu) \,\, dx -0
            \notag \\
    &=\int_{M} \frac{\partial}{\partial t}
            (Fu) \,\, dx
      -\int_{M} \triangle (Fu) \,\, dx
            \notag \\
    &=\int_{M} uF_{t}
       +(Fu_{t}-F\triangle u)
       -u\triangle F
       -2\,\,\nabla u\,\,\nabla F \,\,dx
       \notag\\
    &=\int_{M} -u(\triangle - \partial_{t}) F
       -2\,\,\nabla u\,\,\nabla F \,\,dx  \,.
       \notag
\end{align}
By the result of Lemma \ref{F}, we get,
\begin{align}
    \frac {\partial W}
          {\partial t}
  &=\int_{M} - \frac{a^2\tau}{\pi}u|f_{ij}-\frac{g_{ij}}{2\tau}|^2
    -\big(1-\frac{a^2}{2\pi}\big)u
                      \left(
                            |\nabla f|^2+ \frac{n}{2\tau}
                      \right)
    -\frac{a^2}{\pi}u \triangle f\,dx
\notag\\
&\quad
   -\int_{M}
   \frac{a^2\tau}{\pi}u\, \nabla f \,\nabla (|\nabla f|^2)
   +\frac{a^2\tau}{\pi}u R_{ij}f_{i}f_{j}
   +2\,\,\nabla u\,\,\nabla F \,\,dx
       \notag\\
 &= - \frac{a^2 \tau}{\pi} \int_{M} u
         |f_{ij}-\frac{g_{ij}}{2\tau}|^2 \,\,dx \,\,
    - \big(   1-\frac{a^2}{2\pi}   \big)
    \int_{M} u
                      \left(
                            |\nabla f|^2+ \frac{n}{2\tau}
                      \right)dx
\notag\\
&\quad
   -\int_{M}\frac{a^2}{\pi}u \triangle f
   +\frac{a^2\tau}{\pi}u\, \nabla f \,\nabla (|\nabla f|^2)
   +\frac{a^2\tau}{\pi}u R_{ij}f_{i}f_{j} dx
\notag\\
&\quad
  -\int_{M} 2\,\nabla u\,\frac {a^2\tau} {2\pi}\nabla (|\nabla f|^2)
  +2\,\nabla u\nabla f\,dx\,.
\notag
\end{align}
The last equality comes from $\nabla F=\frac {a^2\tau} {2\pi}\nabla
(|\nabla f|^2) +\nabla f$. Also notice $u \nabla f=-\nabla u $,
there are two more terms canceled. The above becomes,
\begin{align}
 \frac {\partial W}
          {\partial t}
 &= - \frac{a^2 \tau}{\pi} \int_{M} u
         |f_{ij}-\frac{g_{ij}}{2\tau}|^2 \,\,dx \,\,
    -\big(1-\frac{a^2}{2\pi}\big)\int_{M}u
                      \left(
                            |\nabla f|^2+ \frac{n}{2\tau}
                      \right)dx
   \notag
   \\
   &\,\,\,\,\,\,
      +\int_{M}\frac{a^2}{2\pi}u
       \left(
             \frac{2\triangle u}{u}-2|\nabla f|^2
       \right)
     -\frac{a^2\tau}{\pi}u R_{ij}f_{i}f_{j} + 2 u |\nabla f|^2 \,\,dx
       \notag\\
 &= - \frac{a^2 \tau}{\pi} \int_{M} u
         |f_{ij}-\frac{g_{ij}}{2\tau}|^2 \,\,dx \,\,
    +\big(1-\frac{a^2}{2\pi}\big)\int_{M}u
                      \left(
                            -|\nabla f|^2- \frac{n}{2\tau}
                      \right)dx
\notag\\
   &\quad
    -\frac{a^2\tau}{\pi}\int_{M}u R_{ij}f_{i}f_{j}
    +\int_{M}\frac{a^2}{\pi}\triangle u
    + \big(1-\frac{a^2}{2\pi}\big)\int_{M} 2u (|\nabla f|^2)
       \,dx    \,.
\notag
\end{align}
Reorganize the terms, and by the fact $\int_{M} \triangle u\,dx=0$
for closed manifold $M$, we have
\begin{align}
 \frac{\partial W}{\partial t}
 &= - \frac{a^2 \tau}{\pi} \int_{M} u
         |f_{ij}-\frac{g_{ij}}{2\tau}|^2 \,\,dx \,\,
    -\frac{a^2\tau}{\pi}\int_{M}u R_{ij}f_{i}f_{j}\,dx
   \notag
   \\
   &\,\,\,\,\,\,\,
    + \big(1-\frac{a^2}{2\pi}\big)
                     \int_{M} u
                               \left(
                                     |\nabla f|^2
                                     -\frac{\triangle u}{u}
                                     -\frac{n}{2\tau}
                               \right)dx
    +
    \big(1+\frac{a^2}{2\pi}\big)\int_{M} \triangle u \,\,dx
   \notag\\
  &= - \frac{a^2 \tau}{\pi} \int_{M} u
         |f_{ij}-\frac{g_{ij}}{2\tau}|^2 \,\,dx \,\,
    -\frac{a^2\tau}{\pi}\int_{M}u R_{ij}f_{i}f_{j}\,dx
   \notag
   \\
   &\,\,\,\,\,\,\,
    + \big(1-\frac{a^2}{2\pi}\big)
                     \int_{M} u
                               \left(
                                     \frac{|\nabla u|^2}
                                            {u^2}
                                     -\frac{\triangle u}{u}
                                     -\frac{n}{2\tau}
                               \right)dx
   \notag
\end{align}
therefore, for the positive solution $u$ and non-negative Ricci
curvature, all three integral terms are non-positive, the
non-positivity of the last integral term is due to the Li-Yau
\cite{LY} gradient estimate $\frac{|\nabla u|^2} {u^2}
-\frac{\triangle u}{u}-\frac{n}{2\tau}\leq0$, when $\tau=t$.
\end{proof}
Now we give the proof for Theorem \ref{maintheorem-2}, starting with
the following lemma.

\begin{lemma}    \label{P}
Let $u$,\,$f$ be defined as in Lemma \ref{nabla-f}, define
$P:=t(2\triangle f -|\nabla f|^2)u$, then
\begin{equation}   \label{equa-P}
(\triangle - \partial_{t}) P =
   2\,t\,u\,|f_{ij}-\frac{g_{ij}}{2t}|^2
   +u(|\nabla f|^2- \frac{n}{2t})
   +2 \,t\,u\, R_{ij}f_{i}f_{j}  \,.
\end{equation}
\end{lemma}

\begin{proof}
By direct computation, we have
\begin{align}
 (\triangle - \partial_{t}) P
  & =t u\,(\triangle - \partial_{t})(2\triangle f -|\nabla f|^2)
     +t(2\triangle f-|\nabla f|^2)\,\triangle u
\notag\\
  & \quad
     +2 \,t \,\nabla (2\triangle f -|\nabla f|^2)\,\nabla u
     -(2\triangle f -|\nabla f|^2)\, u
     -t \,(2\triangle f -|\nabla f|^2)\,u_{t}\notag  \,.
\end{align}
Since  $\triangle u-u_{t}=0$, two more terms are canceled,
\begin{align}
 (\triangle - \partial_{t}) P
 &=t u\,(\triangle - \partial_{t})(2\triangle f -|\nabla f|^2)
   +2 \,t \,\nabla (2\triangle f -|\nabla f|^2)\,\nabla u
\notag\\
 & \quad
        -(2\triangle f -|\nabla f|^2)\, u  \,.
\notag
\end{align}
Using the result in Lemma \ref{triangle-f}, we have
\begin{align}
 (\triangle - \partial_{t}) P  &=t u\,
  [\,\,
        2|f_{ij}|^2 + 2\,\,\nabla f \,\,\nabla(2\,\triangle f-|\nabla f|^2)
        + 2 R_{ij} f_{i}f_{j}
  \,\,]
\notag\\
  & \quad+2 \,t \,\nabla u\,\nabla (2\triangle f -|\nabla f|^2)
        -(2\triangle f -|\nabla f|^2)\, u \notag   \,.
\end{align}
By the fact $\nabla f=-\frac{\nabla u}{u}$, $\nabla u =-u \,\nabla f
$, the 2nd term and 4th term are canceled, then
\begin{align}
(\triangle - \partial_{t}) P
  & =2\,t u\,|f_{ij}|^2
     +2\,t u\,R_{ij} f_{i}f_{j}
     -(2\triangle f -|\nabla f|^2)\, u   \,.
\end{align}
Keeping in mind that in local normal coordinate $\sum
f_{ij}g_{ij}=\triangle f$, and $\sum g_{ij}^2=n$, by completing the
square,
\begin{align}
(\triangle - \partial_{t}) P
 &=2\,t u\,|f_{ij}-\frac{g_{ij}}{2t}|^2
   -2\,t u\,\left(
                    \frac{g_{ij}}{2t}
                \right)^2
   +2\,t u\,\frac{2f_{ij}g_{ij}}{2t}
\notag\\
&\quad +2\,t u\,R_{ij} f_{i}f_{j}
     -(2\triangle f -|\nabla f|^2)\, u
\notag\\
 &=2\,t u\,|f_{ij}-\frac{g_{ij}}{2t}|^2
   -u\,\frac{n}{2t}+u\,(2\triangle f)
   +2\,t u\,R_{ij} f_{i}f_{j}
   -(2\triangle f -|\nabla f|^2)\, u
\notag\\
&=2\,t u\,|f_{ij}-\frac{g_{ij}}{2t}|^2
   +u(|\nabla f|^2- \frac{n}{2t})
   +2 \,t u\, R_{ij}f_{i}f_{j}
\notag  \,.
\end{align}
\end{proof}
Finally we are in a position to give
\begin{proof}
(of theorem \ref{maintheorem-2}) Define
\begin{equation}
H(f,t,a)
 :=
    \left(
          t\,\big(1+\frac{a^2}{2\pi}\big) \triangle f \,-t\,|\nabla
          f|^2+f-\frac{n}{2}\big(1+\frac{a^2}{2\pi}\big)
    \right)u   \,.
\end{equation}
Reorganize the terms by using $\triangle f =-\frac{\triangle
u}{u}+\frac{|\nabla u|^2}{u^2}=-\frac{\triangle u}{u}+|\nabla
f|^2$,
\begin{align}   \label{reorganize H}
H&=\left(
          \frac{a^2}{2\pi}\,t\, \triangle f
          +t\,(\triangle f-|\nabla f|^2)
          +f-\frac{n}{2}\big(1+\frac{a^2}{2\pi}\big)
    \right)u
 \notag\\
 &=\left(
   \frac{a^2t}{2\pi}(2\triangle f-|\nabla f|^2)
   -\frac{a^2}{2\pi}t (\triangle f-|\nabla f|^2)
+t(\triangle f-|\nabla f|^2)
  \right)u
\\
&\quad+(f-\frac{n}{2})u-\frac{na^2}{4\pi}u \notag  \,.
\end{align}
Combine the last two terms in the brackets,
\begin{align}
H&=\left(
   \frac{a^2t}{2\pi}(2\triangle f-|\nabla f|^2)
   +\big(1-\frac{a^2}{2\pi}\big)t \big(-\frac{\triangle u}{u})
  \right)u+(f-\frac{n}{2})u-\frac{na^2}{4\pi}u
\notag\\
 &=
 \frac{a^2t}{2\pi}(2\triangle f-|\nabla f|^2)u
          -\big(1-\frac{a^2}{2\pi}\big)t \triangle u
            +(f-\frac{n}{2})u-\frac{na^2}{4\pi}u
\notag\\
&=\frac{a^2}{2\pi}P-\big(1-\frac{a^2}{2\pi}\big) t \, \triangle u
            +(f-\frac{n}{2})u-\frac{na^2}{4\pi}u
\notag
\end{align}
where $P=t\,(2\triangle f- |\nabla f|^2) u$ as defined in Lemma
\ref{P}, then,
\begin{align}   \label{heat on H}
(\triangle - \partial_{t}) H
 &=
   \frac{a^2}{2\pi}(\triangle - \partial_{t})P
   -\big(1-\frac{a^2}{2\pi}\big) t
    (\triangle - \partial_{t})\,(\triangle u)
\\
   &\quad
   +\big(1-\frac{a^2}{2\pi}\big)\,\triangle u
   +(\triangle f - f_{t})u+2\,\nabla f \,\nabla u
   \notag   \,.
\end{align}
Notice that
\begin{align}
 \triangle u-u_{t}=0 \,\, &\Rightarrow
 \,\,0=\triangle \big( \triangle u-u_{t} \big)
 =\triangle (\triangle  u)-\triangle (u_{t})
 \notag\\
\,\, &\Rightarrow
 \,\,0=\triangle (\triangle  u)-\partial_{t}(\triangle u)
 =(\triangle -\partial_{t})(\triangle u)
\notag
\end{align}
which means, $\triangle u $ is also a solution to the heat equation;
also observed $\triangle f=f_{t}+|\nabla f|^2+\frac{n}{2t}$ and
$\nabla u=-u\,\nabla f$, then the above expression (\ref{heat on H})
can be simplified as
\begin{align}
(\triangle - \partial_{t}) H &=
   \frac{a^2}{2\pi}(\triangle - \partial_{t})P
   +\big(1-\frac{a^2}{2\pi}\big)\,\triangle u
   +(|\nabla f|^2+\frac{n}{2t})u-2\,|\nabla f|^2\,u
\notag\\
&=
   \frac{a^2}{2\pi}(\triangle - \partial_{t})P
   +\big(1-\frac{a^2}{2\pi}\big)\,\triangle u
   +(\frac{n}{2t}-|\nabla f|^2)u  \,.
\end{align}
By the result in Lemma \ref{P}, we have
\begin{align}
(\triangle - \partial_{t}) H
 &=\left(
   \frac{a^2t}{\pi}|f_{ij}-\frac{g_{ij}}{2t}|^2\,u
   +\frac{a^2}{2\pi}(|\nabla f|^2- \frac{n}{2t})u
   +\frac{a^2t}{\pi}u R_{ij}f_{i}f_{j}
   \right)
\notag\\
 &\quad
   +\big(1-\frac{a^2}{2\pi}\big)\,\frac{\triangle u}{u}\cdot u
   +\big(\frac{n}{2t}-\frac{|\nabla u|^2}{u^2}\big)u
\notag\\
 &=\frac{a^2t}{\pi}|f_{ij}-\frac{g_{ij}}{2t}|^2\,u
   +\frac{a^2t}{\pi}u R_{ij}f_{i}f_{j}
\\
 &\quad
   +\big(\frac{a^2}{2\pi}-1\big)
   \Big(\frac{|\nabla u|^2}{u^2}-\frac{\triangle u}{u}-\frac{n}{2t}\Big)u
   \notag  \,.
\end{align}
If $0 \leq a^2\leq 2\pi$, and $M$ has non-negative Ricci
curvature, then all three terms are non-negative,
\begin{equation}
(\triangle - \partial_{t}) H(f,t,a)\geq 0  \,.
\end{equation}
We claim that,
\begin{equation}  \label{H(0)}
\lim_{t \to 0} H(f,t,a)\leq 0\,.
\end{equation}
Therefore by the maximum principle, for any $t\geq 0 $
\begin{equation}
 H(f,t,a)\leq 0
\notag
\end{equation}
that is,
\begin{equation}
\left(
          t\,\big(1+\frac{a^2}{2\pi}\big) \triangle f \,-t\,|\nabla
          f|^2+f-\frac{n}{2}\big(1+\frac{a^2}{2\pi}\big)
    \right)u \leq 0  \,.
\end{equation}
Recall that $u$ is a positive solution to the heat equation,
consequently
\begin{equation}
    t\,\big(1+\frac{a^2}{2\pi}\big) \triangle f \,-t\,|\nabla
          f|^2+f-\frac{n}{2}\big(1+\frac{a^2}{2\pi}\big)
    \leq 0   \,.
\end{equation}
Let $\alpha=1+\frac{a^2}{2\pi}$, where $0 \leq a^2\leq 2\pi$, then
for any $1\leq \alpha \leq 2$,
\begin{equation}
t\, \left(
    \alpha \,\triangle f \, -|\nabla f|^2
    \right)
+f-\alpha \,\frac{n}{2}\leq 0   \,.
\end{equation}
In particular, if $\alpha=2$, then it becomes the following
differential inequality, which is one of the results proved in
\cite{N1}:
\begin{equation}  \label{Nilei-inequality}
t\,(2 \triangle f \, -|\nabla f|^2)+f-n\leq 0    \,.
\end{equation}
To prove the case for $\alpha >2$, consider Li-Yau gradient
estimate,
\begin{align}
 \frac {|\nabla u|^2}{u^2}-\frac{\triangle u}{u}-\frac{n}{2t}\leq 0
  \,\, &\Rightarrow
 \,\,\triangle f-\frac{n}{2t}\leq 0
  \,\, \Rightarrow\,\,
  t\left(\triangle f-\frac{n}{2t}\right)\leq 0
\notag\\
 &\Rightarrow\,\,t\,\triangle f -\frac{n}{2}\leq0
 \,\, \Rightarrow\,\,(\alpha -2)(t\,\triangle f -\frac{n}{2})\leq0
\notag\\
&\Rightarrow\,\,(\alpha -2)\,t\,\triangle f -\alpha
\,\frac{n}{2}+n\leq0  \label{from-Li-Yau}   \,.
\end{align}
Combine the inequality (\ref{from-Li-Yau}) and one of the results
proved in [N1], that is, the inequality (\ref{Nilei-inequality}), we
obtain
\begin{equation}
   t\,
   \left(\alpha \,\triangle f \, -|\nabla f|^2
   \right)
 +f-\alpha \,\frac{n}{2}\leq 0
 \qquad \text{for any $\alpha > 2$}  \,.
\end{equation}
This proves Theorem \ref{maintheorem-2}, except for the claim
(\ref{H(0)}). The following gives the proof of (\ref{H(0)}):

Let $h(y,t)$ be any positive smooth function with compact support,
by Equation (\ref{reorganize H}), we have
\begin{align}
&\int_{M} h(y,t) \,H(f,t,a) dy
\notag\\
 &=\int_{M}
   \left(
         \frac{a^2t}{2\pi}\Big(2\triangle f-|\nabla f|^2\Big)u
          -\big(1-\frac{a^2}{2\pi}\big) t  \triangle u
            +\big( f-\frac{n}{2} \big)u-\frac{na^2}{4\pi}u
   \right)h dy
\notag\\
 &=  \frac{a^2t}{2\pi}\int_{M} (\triangle f- |\nabla f|^2) u \,h \,dy
     - \big(1-\frac{a^2}{2\pi}\big)t \int_{M} h\,\triangle u\,dy
\notag\\
 & \quad+\frac{a^2t}{2\pi} \int_{M} (\triangle f-\frac{n}{2t})u\,h\, dy
     +\int_{M} (f-\frac{n}{2})u\,h \,dy\,.
\notag
\end{align}
Since $\triangle f=-\frac{\triangle u}{u}+\frac{|\nabla
u|^2}{u^2}$, we have $(\triangle f- |\nabla f|^2) u=\,-\triangle
u$, then
\begin{align}
&\int_{M} h(y,t) \,H(f,t,a) dy
\notag\\
 &=\frac{a^2t}{2\pi}\int_{M} -(\triangle u) \,h \,dy
     + (\frac{a^2}{2\pi}-1)t \int_{M} h\,\triangle u\,dy
\notag\\
 & \quad+ \frac{a^2t}{2\pi}
     \int_{M}
             \left(
              \frac{|\nabla u|^2}{u^2}-\frac{\triangle u}{u} -\frac{n}{2t}
             \right)u \,h \,dy+\int_{M} (f-\frac{n}{2})u\,h\, dy
\notag\\
   &=- t \int_{M} h\,\triangle u\,dy
     + \frac{a^2t}{2\pi}
     \int_{M}
             \left(
               \frac{|\nabla u|^2}{u^2}-\frac{\triangle u}{u} -\frac{n}{2t}
             \right)u\,h dy+\int_{M} (f-\frac{n}{2})u\,h dy
\notag\\
   &= I+II+III\,.  \label{eachterm}
\end{align}
We estimate each term as $t \to 0$,
\begin{align}
 &(I)=- t\int_{M} h\,\triangle u\,dy
     =- t \int_{M} u\,\triangle h\,dy \to 0
        \,\,\,\,as \,\,t \to 0 \,\, (integration\,\, by \,\, parts)
\notag\\
 &(II)=\frac{a^2t}{2\pi}
     \int_{M}
             \left(
               \frac{|\nabla u|^2}{u^2}-\frac{\triangle u}{u} -\frac{n}{2t}
             \right)u \,h \,dy \leq 0  \,\,\,\,as \,\,t
    \to 0
\notag
\end{align}
due to Li-Yau gradient estimate $\frac{|\nabla u|^2}{u^2}
-\frac{\triangle u}{u}-\frac{n}{2t}\leq 0 $.\\
\\
$\noindent\quad(III)\leq 0 \,\,as \,\,t \to 0$ by the following
argument. Let $x$ be a fixed point in $M$, and $(y,t)\in M\times
(0,T]$, by the asymptotic behavior of the fundamental solution
$u=\frac{e^{-f}}{(4\pi\,t)^{\frac{n}{2}}}$ to the heat equation as
$t\to 0$(\cite{MP}),
\begin{align}  \label{asymptotic-u}
 & u(x,y,t)
   \sim
   \frac{e^{-\frac{d^2(x,y)}{4t}}}
        {(4\pi\,t)^{\frac{n}{2}}}
   \sum_{j=0}^{\infty}\tau^{j}u_{j}(x,y,\,t)
\end{align}
where $d(x,y)$ is the distance function. By (\ref{asymptotic-u}) we
mean that there exists a suitable small $T>0$ and a sequence
$(u_{j})_{j\in \mathbf{N}}$ with $u_{j}\in C^{\infty}(M \times M
\times [0,T])$ such that
\begin{align}
 u(x,y,\,t)
 -
 \frac{e^{-\frac{d^2(x,y)}{4\,t}}}
        {(4\pi\,t)^{\frac{n}{2}}}
   \sum_{j=0}^{m}\,t^{j}u_{j}(x,y,\,t)
   :=w_{m}(x,y,\,t)
\end{align}
with
\begin{align}
w_{m}(x,y,\,t)=O(\,t^{m+1-\frac{n}{2}})
\end{align}
as $\,t\to 0$, uniformly for all $x,y\in M$. The function
$u_{0}(x,y,0)$ in (\ref{asymptotic-u}) can be chosen so that
$u_{0}(x,y,0)=1$. Therefore,
\begin{align}
f=
 \frac{d^2(x,y)}{4\,t}
 +
 \ln \left(   1+\,t u_{1}+\,t^2 u_{2}+...+\,t^m u_{m} +O(\,t^{m+1-\frac{n}{2}})  \right)
 \to
 \frac{d^2(x,y)}{4\,t}
 \quad \text{as}\quad \,t\to 0
\end{align}
Thus,
\begin{align}
 &\lim_{\,t \to 0}(III)
 =
 \lim_{\,t \to0}
                \int_{M}(f-\frac{n}{2})u h(y, \,t)\,dy
 =
 \lim_{\,t \to0}
                \int_{M}\left(   \frac{d^2(x,y)}{4\,t}-\frac{n}{2}  \right)u
                h(y, \,t)\,dy
\notag\\
&=
 \lim_{\,t \to0}
                \int_{M}\left(   \frac{d^2(x,y)}{4\,t}-\frac{n}{2}  \right)
                \frac{e^{-\frac{d^2(x,y)}{4\,t}}}
                 {(4\pi\,t)^{\frac{n}{2}}}
                 \left(   1+\,t u_{1}+\,t^2 u_{2}+...+\,t^m u_{m} +O(\,t^{m+1-\frac{n}{2}})  \right)
                h(y, \,t)\,dy
\notag\\
&=
 \lim_{\,t \to0}
                \int_{M}\left(   \frac{d^2(x,y)}{4\,t}-\frac{n}{2}  \right)
                \frac{e^{-\frac{d^2(x,y)}{4\,t}}}
                 {(4\pi\,t)^{\frac{n}{2}}}
                h(y, \,t)\,dy.
\notag
\end{align}

It is easy to see that for any given $\delta>0$, the integration of
the above integrand in the domain $d(x, y) \ge \delta$ converges to
zero exponentially fast. Therefore
\begin{align}
 \lim_{\,t \to 0}(III)
=
 \lim_{\,t \to0}
                \int_{d(x, y) \le \delta}\left(   \frac{d^2(x,y)}{4\,t}-\frac{n}{2}  \right)
                \frac{e^{-\frac{d^2(x,y)}{4\,t}}}
                 {(4\pi\,t)^{\frac{n}{2}}}
                h(y, \,t)\,dy\,.
\end{align}
When $\delta$  is sufficiently small, $d(x, y)$ is sufficiently
close to the Euclidean distance. After a standard approximation
process using local normal coordinates, it is clear that
\begin{align}
 \lim_{\,t \to 0}(III)
=
 \lim_{\,t \to0}
                \int_{R^n}\left(   \frac{|x-y|^2}{4\,t}-\frac{n}{2}  \right)
                \frac{e^{-\frac{|x-y|^2}{4\,t}}}
                 {(4\pi\,t)^{\frac{n}{2}}}
                h_p(y)\,dy\,.
\end{align}
Here $h_p$ is the pull back of $h(\cdot, 0)$ to the Euclidean
space from the region $d(x, y) \le \delta$.

We split the above integral to
\begin{align}
 \lim_{\,t \to 0}(III)
&=
 \lim_{\,t \to0}
                \int_{R^n}\left(   \frac{|x-y|^2}{4\,t}-\frac{n}{2}  \right)
                \frac{e^{-\frac{|x-y|^2}{4\,t}}}
                 {(4\pi\,t)^{\frac{n}{2}}}
                h_p(x)\,dy
\\
&\quad
 +\lim_{\,t \to0}
 \int_{R^n}\left(   \frac{|x-y|^2}{4\,t}-\frac{n}{2}  \right)
                \frac{e^{-\frac{|x-y|^2}{4\,t}}}
                 {(4\pi\,t)^{\frac{n}{2}}}
                (h_p(y)-h_p(x))\,dy.
 \notag
\end{align}
By a straight forward calculation, the second integral on the right
hand side of the last identity converges to zero as $\,t \to 0$,
since $|h_p(y)-h_p(x)| \le C |x-y|$. Hence
\begin{align}
 \lim_{\,t \to 0}(III)
 &= h_p(x)  \lim_{\,t \to0}
                \int_{R^n}\left(   \frac{|x-y|^2}{4\,t}-\frac{n}{2}  \right)
                \frac{e^{-\frac{|x-y|^2}{4\,t}}}
                 {(4\pi\,t)^{\frac{n}{2}}}
                \,dy
\\
 &= h_p(x) \,[\, \lim_{\,t \to0}
                \int_{R^n}   \frac{|y|^2}{4\,t}
                \,
                \frac{e^{-\frac{|y|^2}{4\,t}}}
                 {(4\pi\,t)^{\frac{n}{2}}}
                \,dy - \frac{n}{2} \,]
 =0.
\notag
\end{align}
The last step is by an integration as in an exercise in calculus.

Since all three terms in (\ref{eachterm}) are non-positive as $\,t
\to 0$, we conclude for any positive smooth function $h(y,t)$ with
compact support,
\begin{equation}
  \int_{M} h(y,\,t) \,H(f,\,t,a) d\mu_{\,t}(y)\leq 0
\notag
\end{equation}
then for any small $\,t\geq 0 $
\begin{equation}
 H(f,\,t,a)\leq 0.
\notag
\end{equation}
This completes the proof of the claim \ref{H(0)}.
\end{proof}

\section{Appendix}

The material below, due to Perelman \cite{P},  can be found in
several recent papers and books. They are given for completeness.
Here we follow the presentation in \cite{T}.

\begin{lemma}
Let $u,f,g,\tau$ defined as in Theorem \ref{Ricci}, then the F-
functional defined by
 $F(g,f):=\int \big(R+|\nabla f|^2\big)u\,dx$ is non-decreasing in
 $t$ under
\begin{equation}
\frac{\partial F(g,f)}{\partial t}=2\int |Ric+Hess(f)|^2 u\,dx\,.
\end{equation}
\end{lemma}
\begin{proof}
For simplicity, we denote $\frac{\partial g}{\partial t}=h $,
$\frac{\partial f}{\partial t}=k$. \\
Keeping in mind that $\frac{\partial }{\partial t}g^{ij}=-h^{ij}$,
we may calculate
\begin{equation}
\frac{\partial |\nabla f|^2}{\partial t}
 = -h(\nabla f, \nabla f )
   +2\langle\nabla k, \nabla f\rangle  \,.
\end{equation}
Also
\begin{equation}
\frac{\partial }{\partial t}R
 = -\langle Ric, h\rangle+ \delta^2h-\triangle (tr h)
\end{equation}
where $\delta$ is the divergence operator $\delta:
\Gamma(\otimes^{k}T^{*}M)\rightarrow\Gamma(\otimes^{k-1}T^{*}M)$
defined by: $\delta(T)=-tr_{12}\nabla T$, here $tr_{12}$ means the
trace over the first and second entries of $\nabla T$. Further,
\begin{equation}
\frac{\partial }{\partial t}dx
 =\frac{1}{2}(tr h)dx
\end{equation}
where $dx$ is the volume form
$dx\equiv\sqrt{det(g_{ij})}dx^{1}\wedge ...\wedge dx^{n}$.\\
Therefore,
\begin{align}
\frac{\partial }{\partial t}F(g,f)
 &=\int
     [-h(\nabla f, \nabla f ) +2\langle\nabla k, \nabla f\rangle-\langle Ric, h\rangle+ \delta^2h-\triangle (tr h)]
 udx
\notag\\
&\quad
 +\int(R+|\nabla f|^2)[-k+\frac{1}{2}(tr h)]
 \frac{e^{-f}}{(4\pi\tau)^{\frac{n}{2}}}dx \,.
\end{align}
Three of these terms may be useful addressed by integrating by
parts.\\ First
\begin{equation}
  \int
  2\langle\nabla k, \nabla f\rangle
  \frac{e^{-f}}{(4\pi\tau)^{\frac{n}{2}}} dx
  =\int
      -2k(\triangle f-|\nabla f|^2)
      \frac{e^{-f}}{(4\pi\tau)^{\frac{n}{2}}} dx  \,.
\end{equation}
Second,
\begin{align}
 \int
    (\delta^2h)
     \frac{e^{-f}}{(4\pi\tau)^{\frac{n}{2}}}dx
 &=\int
     \langle
            \delta h, d(\frac{e^{-f}}{(4\pi\tau)^{\frac{n}{2}}})
     \rangle dx
  =\int
     \langle
            h, \nabla d(\frac{e^{-f}}{(4\pi\tau)^{\frac{n}{2}}})
     \rangle dx
\notag\\
 &=\int
  (h(\nabla f, \nabla f)
  -\langle
          Hess(f), h
   \rangle)
  \frac{e^{-f}}{(4\pi\tau)^{\frac{n}{2}}}dx   \,.
\end{align}
Third,
\begin{align}
  \int
     -\triangle (tr h)
      \frac{e^{-f}}{(4\pi\tau)^{\frac{n}{2}}} dx
 &=\int
      -(tr h)\triangle (\frac{e^{-f}}{(4\pi\tau)^{\frac{n}{2}}})dx
\notag\\
 &=\int
      (\triangle f -|\nabla f|^2) (tr h)\frac{e^{-f}}{(4\pi\tau)^{\frac{n}{2}}}dx \,.
\end{align}
Combining these calculations, we find that
\begin{align}
\frac{\partial }{\partial t}F(g,f)
 &=\int
     [-\langle
        Ric, h
       \rangle
      -\langle
        Hess(f),h
       \rangle
       + (\triangle f-|\nabla f|^2)\,(tr h-2 k)
\notag\\
&\qquad\,\,
 +(R+|\nabla f|^2)[-k+\frac{1}{2}(tr h)]
 \,\,\frac{e^{-f}}{(4\pi\tau)^{\frac{n}{2}}}\,dx
\notag\\
 &=\int
     [\langle
        -Ric -Hess(f),h
       \rangle
       + (2\triangle f-|\nabla f|^2+R)(\frac{1}{2}tr h- k)]
   \,u\,dx  \,.  \label{derivative of F}
\end{align}
Notice that under the Ricci flow given by
\begin{align}
\frac{\partial g}{\partial t}&=-2Ric
\\
\frac{\partial f}{\partial t}&=-\triangle f+|\nabla f|^2-R
\end{align}
which is equivalent to another decoupled system of equations by a
proper transform of diffeomorphisms,
\begin{align}
\frac{\partial \hat{g}}{\partial t}
 &=-2\big(Ric(\hat{g})+Hess_{\hat{g}}(\hat{f})\big)
\\
\frac{\partial \hat{f}}{\partial t}
 &=-\triangle_{\hat{g}}f-R_{\hat{g}}
\end{align}
substituting into (\ref{derivative of F}), we have
\begin{align}
\frac{\partial }{\partial t}F(g,f)
 &=2\int
     [\langle
        -Ric -Hess(f), -Ric -Hess(f)
       \rangle
   \,u\,dx
\notag\\
 &=2\int
     |Ric+Hess(f)|^2
   \,u\,dx \,.
\end{align}
\end{proof}

%
\begin{lemma}
Let $g,f,\tau$ defined as in Theorem \ref{Ricci}, define
 $P:=[\tau(2\triangle f-|\nabla f|^2+R)+f-n]u$, then
\begin{equation}
\Box^{*}P=-2\tau|Ric+Hess(f)-\frac{g}{2\tau}|^2u  \,.
\end{equation}
\end{lemma}

\begin{proof}
Consider $P=\frac{P}{u}u$, we find that
\begin{equation}
\Box^{*}P=\Box^{*}(\frac{P}{u}u)
    =\frac{P}{u}\Box^{*}u-u
    \left(
    \frac{\partial}{\partial t}+\triangle
    \right)
    (\frac{P}{u})-2\langle\nabla \frac{P}{u},\nabla u\rangle  \,.
\notag
\end{equation}
Since $\Box^{*}u=0$, and $\nabla f=-\frac{\nabla u}{u}$, we have
\begin{equation}
\frac{\Box^{*}P}{u}
    =-
    \left(
    \frac{\partial}{\partial t}+\triangle
    \right)
    (\frac{P}{u})+2\langle\nabla \frac{P}{u},\nabla f\rangle
\end{equation}
for the first term on the right hand side,
\begin{align}
  - \left(
    \frac{\partial}{\partial t}+\triangle
    \right)
    (\frac{P}{u})
  &=-
    \left(
    \frac{\partial}{\partial t}+\triangle
    \right)
   [\tau(2\triangle f-|\nabla f|^2+R)+f-n]
\notag\\
  &=(2\triangle f-|\nabla f|^2+R)
\notag\\
  &\quad-\tau
    \left(
    \frac{\partial}{\partial t}+\triangle
    \right)
    (2\triangle f-|\nabla f|^2+R)
    -\left(
    \frac{\partial}{\partial t}+\triangle
    \right)f
\notag
\end{align}
using the evolution equation for $f$ in (\ref{evolution-f}) on the
final term, this reduces to
\begin{align}
  - \left(
    \frac{\partial}{\partial t}+\triangle
    \right)
    (\frac{P}{u})
  &=2\triangle f-2|\nabla f|^2+2R-\frac{n}{2\tau}
   \notag\\
  &\quad-\tau
    \left(
    \frac{\partial}{\partial t}+\triangle
    \right)
    (2\triangle f-|\nabla f|^2+R)   \,.
\end{align}
Recall that
\begin{align}    \label{q}
    \left(
    \frac{\partial}{\partial t}+\triangle
    \right)
   & (2\triangle f-|\nabla f|^2+R)
\notag\\
&=4\langle\,Ric,Hess(f)\rangle +\triangle |\nabla f|^2-2Ric(\nabla
f,\nabla f)
\notag\\
 &\quad
  -2\langle\,\nabla f,\nabla(-\triangle f +|\nabla f|^2-R)\rangle
  +2|Ric|^2 \,.
\end{align}
Also
\begin{align}
2\langle\,\nabla \frac{P}{u},\nabla f\rangle =2\tau\langle\,\nabla
(2\triangle f-|\nabla f|^2+R),\nabla f\rangle +2|\nabla f|^2 \,.
\end{align}
Combining the above expressions altogether, we find that
\begin{align}
\frac{\Box^{*}P}{u}
    &=2\triangle f+2R -\frac{n}{2\tau}
     -\tau
     (4\langle\,Ric,Hess(f)\rangle +2|Ric|^2)
\notag\\
    &\quad
    +\tau
         [
         -\triangle |\nabla f|^2+2Ric(\nabla f,\nabla f)+2\langle\,\nabla f,\nabla (\triangle f)\rangle
         ] \,.
\end{align}
The three terms in the square brackets simplified to
$-2|Hess(f)|^2$, so
\begin{align}
\frac{\Box^{*}P}{u}
    &=2\triangle f+2R -\frac{n}{2\tau}
     \notag\\
     &\quad
     -\tau
     [4\langle Ric,Hess(f)\rangle +2|Ric|^2+2|Hess(f)|^2]
\notag\\
    &=2\triangle f+2R -\frac{n}{2\tau}
     -2\tau
     \left(|Ric+Hess(f)|\right)^2
\notag\\
    &=2\langle Ric+Hess(f),g_{ij}\rangle
    -\frac{g_{ij}^2}{2\tau}
     -2\tau
     \left(|Ric+Hess(f)|\right)^2
\notag\\
    &=-2\tau|Ric+Hess(f)-\frac{g}{2\tau}|^2 \,.
\notag
\end{align}

\end{proof}

\emph{Acknowledgement.}  We would like to thank Professor Bennet
Chow and Professor Lei Ni for their support.

\noindent e-mail:  kuang@math.ucr.edu and qizhang@math.ucr.edu
\end{document}